# NONCOMMUTATIVE RELATIVE DE RHAM–WITT COMPLEX VIA THE NORM

ZHOUHANG MAO

Abstract. In [Ill79], Illusie constructed de Rham–Witt complex of smooth $\mathbb{F}_p$-algebras $R$, which computes the crystalline cohomology of $R$, a $\mathbb{Z}_p$-lift of the de Rham cohomology of $R$. There are two different extensions of de Rham–Witt complex: a relative version discovered by Langer–Zink, and a noncommutative version, called *Hochschild–Witt homology*, constructed by Kaledin. The key to Kaledin's construction is his *polynomial Witt vectors*.

In this article, we introduce a common extension of both: relative Hochschild–Witt homology. It is simply defined to be topological Hochschild homology relative to the Tambara functor $\underline{W}(\mathbb{F}_p)$. Adopting Hesselholt's proof of his HKR theorem, we deduce an HKR theorem for relative Hochschild–Witt homology, which relates its homology groups to relative de Rham–Witt complex. We also identify Kaledin's polynomial Witt vectors as the relative Hill–Hopkins–Ravenel norm, which allows us to identify our Hochschild–Witt homology relative to $\mathbb{F}_p$ with Kaledin's Hochschild–Witt homology. As a consequence, we deduce a comparison between Hochschild–Witt homology and topological restriction homology, fulfilling a missing part of [Kal19].

## 1. Introduction

Let $k$ be a perfect $\mathbb{F}_p$-algebra, and $A$ a smooth $k$-algebra. Suppose that the smooth $k$-algebra $A$ admits a smooth lift to the ring $W(k)$ of Witt vectors[1.1]. In this case, Grothendieck observed that the $p$-completed de Rham cohomology of $\tilde{A}$ over $W(k)$ does not depend on the choice of smooth lift $\tilde{A}$. Berthelot, following Grothendieck's ideas in [GGK+68], gives an intrinsic (i.e. without choice of $\tilde{A}$), site-theoretic description of the $p$-completed de Rham cohomology of $\tilde{A}$ over $W(k)$, called the *crystalline cohomology*, in [Ber74]. In [Ill79], Illusie constructed an $F$-$V$-pro-complex $(W.\Omega_A^*, d)$, called the (absolute) de Rham–Witt complex[1.2], which roughly consists of the following data:

1. for every $r \in \mathbb{N}_{>0}$, a cochain complex $(W_r \Omega_A^*, d)$ of $W_r(k)$-modules;
2. for every $r \in \mathbb{N}_{>0}$, a map $R \colon (W_{r+1} \Omega_A^*, d) \to (W_r \Omega_A^*, d)$, called the *restriction*, of cochain complexes of $W_{r+1}(k)$-modules;
3. for every $r \in \mathbb{N}_{>0}$, two maps $F \colon W_{r+1} \Omega_A^* \to W_r \Omega_A^*$ and $V \colon W_r \Omega_A^* \to W_{r+1} \Omega_A^*$, called the Frobenius map and the Verschiebung map respectively, of graded abelian groups;
4. for every $r \in \mathbb{N}_{>0}$, a map $\lambda_r \colon W_r(A) \to W_r \Omega_A^0$ of $W_r(k)$-algebras,

which satisfy certain conditions. When $r = 1$, the cochain complex $(W_1 \Omega_A^*, d)$ coincides with the de Rham complex $(\Omega_A^*, d)$, thus the de Rham–Witt complex can be viewed as a thickening of the de Rham complex. The key property is that, for every $r \in \mathbb{N}_{>0}$, the complex $(W_r \Omega_A^*, d)$ represents the crystalline cohomology of $A$ over

---

[1.1]. Such a lift always exists, by a theorem of Arabia, cf. [Sta21, Tag 07M8]. However, such a existence is not essential to our discussion.

[1.2]. It is absolute in the sense that it does not quite depend on the base perfect $\mathbb{F}_p$-algebra $k$.

$W_r(k)$, and the map $R\colon(W_{r+1}\Omega_A^*,\mathrm{d})\to(W_r\Omega_A^*,\mathrm{d})$ represents the canonical reduction map between crystalline cohomologies of $A$ induced by the map $R\colon W_{r+1}(k)\to W_r(k)$ of PD-thickenings of $k$.

There are two different extensions of this picture:

1. Langer–Zink extended de Rham–Witt complex to mixed characteristic relative situation in [LZ04]. More precisely, for every $p$-completely smooth map $k\to A$ of ($p$-complete) commutative algebras, they defined the relative de Rham–Witt complex $(W.\Omega_{A/k}^*,\mathrm{d})$ with a similar $F$-$V$-pro-complex structure as above.

2. Let $k$ be a commutative ring, and $A$ a smooth $k$-algebra. Then the Hochschild–Kostant–Rosenberg theorem tells us that, for every $n\in\mathbb{Z}$, the $n$-th Hochschild homology $\mathrm{HH}_n(A/k)$ coincides with the module $\Omega_{A/k}^n$ of differential $n$-forms. Moreover, for every $n\in\mathbb{Z}$, the map $\mathrm{HH}_n(A/k)\to\mathrm{HH}_{n+1}(A/k)$ induced by the $\mathbb{T}$-action on the Hochschild homology $\mathrm{HH}(A/k)$ coincides with the de Rham differential $\Omega_{A/k}^n\to\Omega_{A/k}^{n+1}$. Since the Hochschild homology is well defined for every associative $k$-algebra, we can understand Hochschild homology as a noncommutative version of de Rham complex.

   In [Kal19], Kaledin found a noncommutative version of (absolute) de Rham–Witt complex, called *Hochschild–Witt homology*. More precisely, let $A$ be an associative $\mathbb{F}_p$-algebra. For every $r\in\mathbb{N}_{>0}$, there is a $W(k)$-module spectrum $W_r\mathrm{HH}(A)$ with $\mathbb{T}$-action, such that, letting $k$ be a perfect $\mathbb{F}_p$-algebra and $A$ a smooth $k$-algebra, for every $r\in\mathbb{N}_{>0}$ and every $n\in\mathbb{Z}$, the $n$-th $r$-truncated Hochschild–Witt homology $\pi_n(W_r\mathrm{HH}(A))$ is isomorphic to the $W_r(k)$-module $W_r\Omega_A^n$. Moreover, the map $\pi_n(W_r\mathrm{HH}(A))\to\pi_{n+1}(W_r\mathrm{HH}(A))$ induced by the $\mathbb{T}$-action on $W_r\mathrm{HH}(A)$ coincides with the differential $W_r\Omega_A^n\to W_{r+1}\Omega_A^n$. There are also analogues maps $R$, Frobenius and Verschiebung on Hochschild–Witt homology.

At this point, there is already a natural question: is there a noncommutatve version of relative de Rham–Witt complex à la Langer–Zink? The goal of this article is to give and explain an affirmative answer to this question:

**Definition 1.1.** (Definition 4.1) *Let $k$ be a commutative ring, and $A$ an associative $k$-algebra. Then the* Hochschild–Witt homology $\underline{W}\mathrm{HH}(A/k)$ *of $A$ relative to $k$ is defined to be the topological Hochschild homology*

$$\mathrm{THH}(A)\otimes_{\mathrm{THH}(k)}^{\mathbb{L}}\underline{W}(k)$$

*relative to the pre-$p$-cyclotomic base $\underline{W}(k)$, and $W_r\mathrm{HH}(A/k):=\underline{W}\mathrm{HH}(A/k)^{C_{p^{r-1}}}$.*

We now explain why this coincides with Kaledin's notion when $k=\mathbb{F}_p$. We first summarize Kaledin's construction briefly. The key to Kaledin's Hochschild–Witt complex is his *polynomial Witt vectors*. Let $\tilde{M}$ be a free abelian group, and let $M$ denote its base change $\tilde{M}\otimes_{\mathbb{Z}}\mathbb{F}_p$ along the quotient map $\mathbb{Z}\to\mathbb{F}_p$. Then the module $W_r(\mathbb{F}_p;M)$ of $r$-truncated polynomial Witt vectors is given by the $C_{p^r}$-Tate cohomology $\hat{H}^0(C_{p^r};\tilde{M}^{\otimes p^{r-1}})$, where the $C_{p^r}$-action on $M^{\otimes p^{r-1}}$ is inflated from the $C_{p^{r-1}}$-action of permutation of tensor factors. Kaledin proved that this construction does not depend on the choice of free $\mathbb{Z}$-lift $\tilde{M}$. Moreover, there is a trace theory structure on polynomial Witt vector functors, which allows us to talk about twisted Hochschild homology with respect to them. They give rise to (truncated) Hochschild–Witt homology.

Our first observation is that, similarly to crystalline cohomology, the polynomial Witt vector functor can be intrinsically described without free $\mathbb{Z}$-lift $\tilde{M}$:

**Proposition 1.2. (Corollary 6.6)** *Let $M$ be an $\mathbb{F}_p$-vector space, and $r \in \mathbb{N}_{>0}$. Then the group of $r$-truncated polynomial Witt vectors can be functorially identified with the $C_{p^{r-1}}$-fixed points*

$$\left(M^{\otimes^{\mathbb{L}}_{\underline{W}(\mathbb{F}_p)} C_{p^{r-1}}}\right)^{C_{p^{r-1}}}$$

*of the (Hill–Hopkins–Ravenel) $C_{p^{r-1}}$-norm $M^{\otimes^{\mathbb{L}}_{\underline{W}(\mathbb{F}_p)} C_{p^{r-1}}}$ relative to the Tambara functor $\underline{W}(\mathbb{F}_p)$.*

The key to this identification is a fairly simple description of the norm relative to the constant Tambara functor $\underline{\mathbb{Z}}$:

**Proposition 1.3. (Proposition 5.6 and Lemma 5.11)** *Let $M$ be a free abelian group, and $G$ a finite group. Then the $G$-norm $M^{\otimes^{\mathbb{L}}_{\underline{\mathbb{Z}}} G}$ relative to the Tambara functor $\underline{\mathbb{Z}}$ coincides with the cohomological Mackey functor corresponding to the $G$-module $M^{\otimes_{\mathbb{Z}} G}$.*

Moreover, modifying Kaledin's "first" nontrivial example of trace theory, we show that every Tambara functor gives rise to a trace theory in Section 7. Combining these two, we see that Kaledin's Hochschild–Witt homology can be easily read from our relative Hochschild–Witt homology. More precisely, let $A$ be an associative $\mathbb{F}_p$-algebra, then the $r$-truncated Hochschild–Witt homology $W_r \operatorname{HH}(A)$ is equivalent to the genuine $C_{p^{r-1}}$-fixed points

$$\left(\operatorname{THH}(A) \otimes^{\mathbb{L}}_{\operatorname{THH}(\mathbb{F}_p)} \underline{W}(\mathbb{F}_p)\right)^{C_{p^{r-1}}}$$

of the Hochschild–Witt homology $\underline{W} \operatorname{HH}(A/\mathbb{F}_p)$ of $A$ relative to $\mathbb{F}_p$.

This description also makes the Frobenius map and the Verschiebung map transparent: they simply come from the restriction map and the transfer map of the genuine equivariant homotopy structure. The restriction map on the Hochschild–Witt homology, on the other hand, comes from the pre-$p$-cyclotomic structure on the Hochschild–Witt homology.

This is indeed a noncommutative version of relative de Rham–Witt homology, by the following HKR-type theorem.

**Proposition 1.4. (Proposition 4.9)** *Let $k$ be a commutative ring, and $A$ a smooth commutative $k$-algebra. Then we have isomorphism*

$$\pi_n W_r \operatorname{HH}(A/k) \cong W_r \Omega^n_{A/k}$$

*which is compatible with the Verschiebung $V$, the Frobenius $F$, the restriction $R$, and the $\mathbb{T}/C_{p^{r-1}}$-action on the left hand side corresponds to the differential on the right hand side.*

We basically follow the proof of the HKR theorem in [Hes96] (replacing facts in char $p$ by their counterparts in mixed characteristic): first to produce a structure of $F$-$V$-procomplex on the left hand side, then the universal property of relative de Rham–Witt complex gives rise to a map from the right hand side to the left hand side. Hesselholt's original proof shows that this map is an equivalence on polynomial algebras, and then conclude by the fact that both sides satisfy étale descent and étale base change.

Finally, the direct occurrence of topological Hochschild homology in this description makes it much easier to establish an equivalence between topological restriction homology and (non-truncated) Hochschild–Witt homology for associative $\mathbb{F}_p$-algebras:

**Proposition 1.5. (Corollary 4.2)** *Let $A$ be an associative $\mathbb{F}_p$-algebra. Then there is a canonical equivalence*

$$\mathrm{TR}(A) \longrightarrow \lim_{r,R} W_r \mathrm{HH}(A/\mathbb{F}_p)$$

*of topological Cartier modules.*

This addresses one of missing comparisons mentioned at the end of [Kal19, §0].

**Notation.** *Let $G$ be a finite group. We will denote by $\mathrm{Sp}^{gBG}$ the $\infty$-category of $G$-spectra[1.3], by $\mathrm{Sp}^{g_p B\mathbb{T}}$ the $\infty$-category of $p$-cyclonic spectra (defined in [BG16]).*

*Acknowledgments.* We would like to thank Wolfgang STEIMLE for informing us the Hill–Hopkins–Ravenel norm, and Kaif HILMAN for explaining us the relative norm, without which this article cannot appear. We would also like to thank Victor SAUNIER for explaining trace theory, which is also essential in this article. In addition, we would like to thank Lukas BRANTNER, Bastiaan CNOSSEN, Marc HOYOIS, Guchuan LI, Sil LINSKENS, Jonas MCCANDLESS, Maxime RAMZI, and Mingcong ZENG. This project has received funding from the European Research Council (ERC) under the European Union's Horizon 2020 research and innovation programme (grant agreement No. 864145).

## 2. Relative $p$-cyclotomic spectra

In [BMY23], the authors introduce the concept of pre-$p$-cyclotomic base, which should be a correct home to talk about relative cyclotomic structures. In this section, we mainly discuss a weaker version of this, namely, that of pre-$p$-cyclotomic rings, which simply omits the norm structure on pre-$p$-cyclotomic bases. This concept is enough for us to talk about relative cyclotomic and polygonic modules.

**Definition 2.1.** (([BM15, §4] & [BG16, Def 3.20]) )

1. A **pre-$p$-cyclotomic spectrum** *(or a $p$-typical pre-cyclotomic spectrum)* is a $p$-cyclonic spectrum $X$ equipped with a map $X^{\Phi C_p} \to X$ of $p$-cyclonic spectra. The $\infty$-category $\mathrm{CycSp}_p^{\mathrm{pre}}$ of pre-$p$-cyclotomic spectra *(or)* is the $\infty$-category $\mathrm{Alg}_{(-)^{\Phi C_p}}(\mathrm{Sp}^{g_p B\mathbb{T}})$ of $(-)^{\Phi C_p}$-algebras in $p$-cyclonic spectra. Commutative ring objects there are called **pre-$p$-cyclotomic rings**.

   Since the endofunctor $(-)^{\Phi C_p}$ carries a symmetric monoidal structure, the $\infty$-category $\mathrm{CycSp}_p^{\mathrm{pre}}$ inherits a symmetric monoidal structure from $\mathrm{Sp}^{g_p B\mathbb{T}}$.

2. A **pre-$p$-polygonic spectrum** *(or a $p$-typical pre-polygonic spectrum)* is a sequence $(X_r \in \mathrm{Sp}^{gBC_{p^r}})_{r \in \mathbb{N}}$ of genuine $C_{p^*}$-spectra along with transition maps $X_{r+1}^{\Phi C_p} \to X_r$ for every $r \in \mathbb{N}$. The symmetric monoidal $\infty$-category $\mathrm{PgcSp}_p^{\mathrm{pre}}$ of pre-$p$-polytonic spectra is the lax sequential limit of the tower

$$\cdots \xrightarrow{(-)^{\Phi C_p}} \mathrm{Sp}^{gBC_{p^2}} \xrightarrow{(-)^{\Phi C_p}} \mathrm{Sp}^{gBC_p} \xrightarrow{(-)^{\Phi C_p}} \mathrm{Sp}$$

---

[1.3]. This is usually denoted by $\mathrm{Sp}_G$, $\mathrm{Sp}^G$ or $\mathrm{Sp}^{gG}$ in the literature, but we need the concept of genuine $X$-spectra for every groupoid $X$, thus we opt for the cumbersome notation $\mathrm{Sp}^{gBG}$ in this article.

of symmetric monoidal $\infty$-categories. Commutative ring objects there are called pre-$p$-polygonic rings.

3. Let $A$ be a pre-$p$-cyclotomic ring. We say that an $A$-module $M$ in pre-$p$-cyclotomic spectra is $A$-cyclotomic if the canonical map

$$M^{\Phi C_p} \otimes^{\mathbb{L}}_{A^{\Phi C_p}} A \longrightarrow M$$

of $p$-cyclonic spectra is an equivalence. We will denote by $\mathrm{CycSp}_{p,A}$ the $\infty$-category of $p$-typical $A$-cyclotomic spectra.

4. Let $A$ be a pre-$p$-polygonic ring. We say that an $A$-module $M$ in pre-$p$-polygonic spectra is $A$-polygonic if, for every $r \in \mathbb{N}$, the canonical map

$$M_{r+1}^{\Phi C_p} \otimes^{\mathbb{L}}_{A_{r+1}^{\Phi C_p}} A_r \longrightarrow M_r$$

of $C_{p^r}$-spectra is an equivalence. We will denote by $\mathrm{PgcSp}_{p,A}$ the $\infty$-category of $p$-typical $A$-polygonic spectra.

*Remark* 2.2. There is a forgetful functor $\mathrm{CycSp}_p^{\mathrm{pre}} \to \mathrm{PgcSp}_p^{\mathrm{pre}}, X \mapsto (X_r = X)_{r \in \mathbb{N}}$ which carries a canonical symmetric monoidal structure. Thus a pre-$p$-cyclotomic ring gives rise to a pre-$p$-polygonic ring.

*Remark* 2.3. Unravelling definitions (and using the symmetric monoidal structure on $(-)^{\Phi C_p}$), a pre-$p$-cyclotomic ring can be alternatively given by a commutative algebra in $p$-cyclonic spectra $R$ equipped with a map $R^{\Phi C_p} \to R$ of commutative algebras in $p$-cyclonic spectra. In a private conversation, Allen YUAN told us that, the model-independent language, a pre-$p$-cyclomic base is given by a $\mathbb{T}$-$\mathbb{E}_\infty$-ring $R$ equipped with a map $R^{\Phi C_p} \to R$ of $\mathbb{T}$-$\mathbb{E}_\infty$-rings.

*Remark* 2.4. ([McC21, **Ex 3.1.12**]) Recall that there is an adjunction

$$\mathrm{Sp}^{g_p B\mathbb{T}} \xrightleftharpoons[\mathrm{Infl}_p]{(-)^{\Phi C_p}} \mathrm{Sp}^{g_p B(\mathbb{T}/C_p)} \tag{2.1}$$

where the left adjoint $(-)^{\Phi C_p}$ is symmetric monoidal. Thus the symmetric monoidal $\infty$-category $\mathrm{CycSp}_p^{\mathrm{pre}}$ can be alternatively described as the symmetric monoidal $\infty$-category $\mathrm{CoAlg}_{\mathrm{Infl}_p}(\mathrm{Sp}^{g_p B\mathbb{T}})$ of $\mathrm{Infl}_p$-coalgebras in $p$-cyclonic spectra, i.e. a $p$-cyclonic spectrum $X \in \mathrm{Sp}^{g_p B\mathbb{T}}$ equipped with a map $X \to \mathrm{Infl}_p X$ of $p$-cyclonic spectra. Since the functor $\mathrm{Infl}_p$ is $t$-exact, it follows that the $\infty$-category $\mathrm{CycSp}_p^{\mathrm{pre}}$ inherits a $t$-structure from that on $\mathrm{Sp}^{g_p B\mathbb{T}}$, such that the forgetful functor $\mathrm{CycSp}_p^{\mathrm{pre}} \to \mathrm{Sp}^{g_p B\mathbb{T}}$ is $t$-exact. Moreover, the fully faithful functor $\mathrm{CycSp}_p \to \mathrm{CycSp}_p^{\mathrm{pre}}$ carries the connective part to the connective part.

Similarly, a pre-$p$-polygonic spectrum can be alternatively described as a sequence $(X_r \in \mathrm{Sp}^{gBC_{p^r}})_{r \in \mathbb{N}}$ of genuine $C_{p^*}$-spectra along with transition maps $X_{r+1} \to \mathrm{Infl}_p X_r$ for every $r \in \mathbb{N}$, and it carries a $t$-structure such that the forgetful functor $\mathrm{PgcSp}_p^{\mathrm{pre}} \to \mathrm{Sp}^{gBC_{p^r}}$ is $t$-exact for every $r \in \mathbb{N}$.

*Example* 2.5. By definition, every $p$-cyclotomic spectrum is a pre-$p$-cyclotomic spectrum, and every $p$-cyclotomic ring is a pre-$p$-cyclotomic ring.

The pre-$p$-cyclotomic structure is enough to talk about its topological restriction homology (TR).

*Construction* 2.6. (**[ABG+18, §3.2]**) Let $X$ be a pre-$p$-cyclotomic spectrum. Then for every $r \in \mathbb{N}_{>0}$, we have a composite map
$$R \colon X^{C_{p^{r+1}}} = (X^{C_p})^{C_{p^r}} \longrightarrow (X^{\Phi C_p})^{C_{p^r}} \longrightarrow X^{C_{p^r}}$$
of $p$-cyclonic spectra, which gives rise to a tower
$$\cdots \xrightarrow{R} X^{C_{p^r}} \xrightarrow{R} \cdots \xrightarrow{R} X^{C_p} \xrightarrow{R} X$$
of $p$-cyclonic spectra. The topological restriction homology $\mathrm{TR}(X)$ is defined to be the limit of this tower. Note that, by construction, there is a canonical equivalence $\mathrm{TR}(X)^{C_p} \simeq \mathrm{TR}(X)$, which gives rise to a topological $p$-Cartier module structure on $\mathrm{TR}(X)$. This gives rise to a functor $\mathrm{TR} \colon \mathrm{CycSp}_p^{\mathrm{pre}} \to \mathrm{TCart}_p$.

More formally, let $\mathcal{C}$ be an $\infty$-category with sequential limits, and $F \Rightarrow G$ a natural transformation between two endofunctors $F$ and $G$ of $\mathcal{C}$. Then we have a functor $\mathrm{Alg}_G(\mathcal{C}) \to \mathrm{Fix}_F(\mathcal{C})$ concretely given by
$$X \longmapsto \lim(\cdots \to FX \to X),$$
where the map $FX \to X$ is the composite map $FX \to GX \to X$ (cf. [AN21, Prop 5.2]). We apply this recipe to the case that $\mathcal{C} = \mathrm{Sp}^{g_p B\mathbb{T}}$, $F = (-)^{C_p}$, and $G = (-)^{\Phi C_p}$.

*Construction* 2.7. (**Cartier modules**) Recall that a *$p$-Cartier module* is an abelian group $M$ equipped with two endomorphisms $(V, F) \in \mathrm{End}_{\mathrm{Ab}}(M)^2$ of $M$, called the *Verschiebung* and the *Frobenius*, such that $FV = p$. Let $\mathrm{Cart}_p^{\flat}$ denote the category of $p$-Cartier modules with injective Verschiebung. We have a composite functor
$$\mathrm{Cart}_p^{\flat} \hookrightarrow \mathrm{TCart}_p \longrightarrow \mathrm{CycSp}_p^{\mathrm{gen}} \hookrightarrow \mathrm{CycSp}_p^{\mathrm{pre}} \xrightarrow{\pi_0} \mathrm{CycSp}_p^{\mathrm{pre},\heartsuit},$$
where $\mathrm{TCart}_p$ is the $\infty$-category of *topological $p$-Cartier modules* [AN21]. We denote by $\underline{M}$ the pre-$p$-cyclotomic spectrum associated to the $p$-Cartier module $M$. Unravelling definitions, we get $\underline{M}^{C_{p^{r-1}}} := \mathrm{coker}\left(M \xrightarrow{V^r} M\right) =: M/V^r$, with the Lewis diagram of $\underline{M} \to \mathrm{Infl}_p \underline{M}$ depicted as

$$\begin{array}{ccc}
M/V^{r+1} & \xrightarrow{\mathrm{can}} & M/V^r \\
{\scriptstyle F}\downarrow\uparrow{\scriptstyle V} & & {\scriptstyle F}\downarrow\uparrow{\scriptstyle V} \\
M/V^r & \xrightarrow{\mathrm{can}} & M/V^{r-1} \\
\vdots & & \vdots \\
M/V^2 & \xrightarrow{\mathrm{can}} & M/V \\
{\scriptstyle F}\downarrow\uparrow{\scriptstyle V} & & {\scriptstyle 0}\downarrow\uparrow{\scriptstyle 0} \\
M/V & \xrightarrow{0} & 0
\end{array}.$$

This functor carries a symmetric monoidal structure, where $\mathrm{Cart}_p^{\flat}$ is equipped with the box product, and in particular, a $p$-Cartier commutative algebra with injective Verschiebung[2.1] will give rise to a pre-$p$-cyclotomic ring.

---

[2.1]. This construction works without injectivity of Verschiebung, but the result is morally correct only under this assumption.

*Example* 2.8. **(Witt vectors)** Let $A$ be an associative ring. Recall that $p$-typical Witt vectors $W(A)$ of $A$ form a $p$-Cartier module, and as such, it is isomorphic to $\pi_0 \operatorname{TR}(A)$ (established in [Hes97], cf. [AN21, Thm 6.1]). Let $\underline{W}(A)$ denote the pre-$p$-cyclotomic spectrum corresponding to the $p$-Cartier module $W(A)$ obtained by Construction 2.7, Concretely, we have $\underline{W}(A)^{C_{p^{r-1}}} \cong W_r(A)$, with the Lewis diagram of $\underline{W}(A) \to \operatorname{Infl}_p \underline{W}(A)$ depicted as

$$\begin{array}{ccc}
W_{r+1}(A) & \xrightarrow{R} & W_r(A) \\
F \downarrow \uparrow V & & F \downarrow \uparrow V \\
W_r(A) & \xrightarrow{R} & W_{r-1}(A) \\
\vdots & & \vdots \\
W_2(A) & \xrightarrow{R} & A \\
F \downarrow \uparrow V & & 0 \downarrow \uparrow 0 \\
A & \xrightarrow{0} & 0
\end{array}.$$

Since the inclusion $\operatorname{CycSp}_p \to \operatorname{CycSp}_p^{\mathrm{pre}}$ preserves 1-connective part, we have $\underline{W}(A) \simeq \pi_0(\operatorname{THH}(A))$ as pre-$p$-cyclotomic spectra. Moreover, the association $\operatorname{Alg}(\operatorname{Ab}) \to \operatorname{CycSp}_p^{\mathrm{pre},\heartsuit}$ is symmetric monoidal, thus carries a commutative ring to a pre-$p$-cyclotomic ring. In particular, the pre-$p$-cyclotomic spectrum $\underline{W}(k)$ is canonocially a pre-$p$-cyclotomic ring. In fact, this argument also shows that $\underline{W}(k)$ has a canonical pre-$p$-cyclotomic base structure (as defined in Remark 2.3).

*Example* 2.9. (**[BMY23, §7.3]**) There is a canonical pre-$p$-cyclotomic ring structure on $\mathbb{F}_p$. We give a slightly different aspect to [BMY23, §7.3]. Indeed, let $k$ be a perfect $\overline{\mathbb{F}}_p$-algebra. As a special case of Example 2.8, the $\mathbb{T}$-Mackey functor $\underline{W}(k)$ is canonically a pre-$p$-cyclotomic ring, thus so is its modulo $p$ reduction $\underline{W}(k)/^{\mathbb{L}}p$ and its zero homotopy group $\pi_0(\underline{W}(k)/^{\mathbb{L}}p)$, which is precisely $\underline{k}$. In particular, this gives rise to a map $\underline{W}(k) \to \underline{k}$ of pre-$p$-cyclotomic rings (and as in Example 2.8, it is also a map of pre-$p$-cyclotomic bases).

Unwinding definitions, and the symmetric monoidal structure on $(-)^{\Phi C_p}$, we get

**Lemma 2.10.** *Let $A \to B$ be a map of pre-$p$-cyclotomic (resp. pre-$p$-polygonic) rings. Then the base change functor*

$$\operatorname{Mod}_A(\operatorname{CycSp}_p^{\mathrm{pre}}) \xrightarrow{(-) \otimes_A^{\mathbb{L}} B} \operatorname{Mod}_B(\operatorname{CycSp}_p^{\mathrm{pre}})$$

*(resp.*

$$\operatorname{Mod}_A(\operatorname{PgcSp}_p^{\mathrm{pre}}) \xrightarrow{(-) \otimes_A^{\mathbb{L}} B} \operatorname{Mod}_B(\operatorname{PgcSp}_p^{\mathrm{pre}})$$

*) restricts to a functor $\operatorname{CycSp}_{p,A} \to \operatorname{CycSp}_{p,B}$ (resp. $\operatorname{PgcSp}_{p,A} \to \operatorname{PgcSp}_{p,B}$).*

## 3. $p$-CYCLOTOMIC SPECTRA RELATIVE TO $\underline{W}(\mathbb{F}_p)$

In this section, we will show that the TR of pre-$p$-cyclotomic spectra is invariant under base change along $\operatorname{THH}(\mathbb{F}_p) \to \underline{W}(\mathbb{F}_p)$, which is the key to compare our Hochschild–Witt homology relative to $\mathbb{F}_p$ with TR in Section 4.

**Notation 3.1.** *Let $M$ be a $p$-cyclonic spectrum. Then we denote by $\operatorname{TF}(M)$ the sequential limit of the tower*

$$\cdots \longrightarrow M^{C_{p^2}} \longrightarrow M^{C_p} \longrightarrow M$$

*along the restriction maps*[3.1]. *In particular, for any ring spectrum $R$, we will denote by* $\mathrm{TF}(R)$ *the spectrum* $\mathrm{TF}(\mathrm{THH}(R))$.

**Proposition 3.2.** *Let $M$ be a $\mathrm{THH}(\mathbb{F}_p)$-module in pre-p-cyclotomic spectra (resp. pre-p-polygonic spectra). Then the canonical map*

$$M \longrightarrow M \otimes^{\mathbb{L}}_{\mathrm{THH}(\mathbb{F}_p)} \underline{W}(\mathbb{F}_p)$$

*of pre-p-cyclotomic spectra becomes an equivalence of topological Cartier modules (resp. spectra) after applying* $\mathrm{TR}$. *In particular, we have an equivalence* $\mathrm{TR}(\underline{W}(\mathbb{F}_p)) \simeq \mathrm{TR}(\mathbb{F}_p)$ *of topological Cartier modules.*

*Proof.* We prove the statement in the cyclotomic case. The polygonic case is similar. By [HM97] (cf. [Sul20, Prop 4.7] for a perfectoid generalization), the Bökstedt element $\sigma \in \pi_2 \mathrm{THH}(\mathbb{F}_p)$ lifts to a Bökstedt element $\sigma \in \pi_2 \mathrm{TF}(\mathbb{F}_p)$, and the map on homotopy groups induced by the restriction map

$$R \colon \mathrm{TF}(\mathbb{F}_p) \longrightarrow \mathrm{TF}(\mathbb{F}_p)$$

is given by the map

$$\begin{aligned}\mathbb{Z}_p[\sigma] &\longrightarrow \mathbb{Z}_p[\sigma] \\ \sigma &\longmapsto p\lambda\sigma\end{aligned}$$

of graded $\mathbb{Z}_p$-algebras for some unit $\lambda \in \mathbb{Z}_p^\times$. It follows that, the restriction maps fit into a commutative diagram

$$\begin{array}{ccc} \mathrm{THH}(\mathbb{F}_p)^{C_p}[2] & \xrightarrow{\sigma} & \mathrm{THH}^{C_p}(\mathbb{F}_p) \\ \downarrow{p\lambda R} & & \downarrow{R} \\ \mathrm{THH}(\mathbb{F}_p)[2] & \xrightarrow{\sigma} & \mathrm{THH}(\mathbb{F}_p) \end{array}$$

in $\mathrm{Sp}^{g_p B\mathbb{T}}$. It follows that the Postnikov truncation $\mathrm{THH}(\mathbb{F}_p) \to \tau_{\leqslant 0}(\mathrm{THH}(\mathbb{F}_p)) = \underline{W}(\mathbb{F}_p)$ in $\mathrm{CAlg}(\mathrm{CycSp}_p^{\mathrm{pre}})$ fits into a fiber sequence

$$\mathrm{THH}(\mathbb{F}_p)[2] \xrightarrow{\sigma} \mathrm{THH}(\mathbb{F}_p) \longrightarrow \underline{W}(\mathbb{F}_p) \qquad (3.1)$$

in $\mathrm{Mod}_{\mathrm{THH}(\mathbb{F}_p)}(\mathrm{CycSp}_p^{\mathrm{pre}})$, where the pre-$p$-cyclotomic structure on $\mathrm{THH}(\mathbb{F}_p)[2]$ is given by the composite map

$$\mathrm{THH}(\mathbb{F}_p)^{\Phi C_p}[2] \xrightarrow{\simeq} \mathrm{THH}(\mathbb{F}_p)[2] \xrightarrow{p\lambda} \mathrm{THH}(\mathbb{F}_p)[2] \qquad (3.2)$$

in $\mathrm{Mod}_{\mathrm{THH}(\mathbb{F}_p)}(\mathrm{Sp}^{g_p B\mathbb{T}})$. We note that, the tower

$$\cdots \to (M \otimes^{\mathbb{L}}_{\mathrm{THH}(\mathbb{F}_p)} \mathrm{THH}(\mathbb{F}_p)[2])^{C_{p^{r+1}}} \to (M \otimes^{\mathbb{L}}_{\mathrm{THH}(\mathbb{F}_p)} \mathrm{THH}(\mathbb{F}_p)[2])^{C_{p^r}} \to \cdots \qquad (3.3)$$

of spectra computing $\mathrm{TR}(M \otimes^{\mathbb{L}}_{\mathrm{THH}(\mathbb{F}_p)} \mathrm{THH}(\mathbb{F}_p)[2])$ is pro-truncated zero, since $(N \otimes^{\mathbb{L}}_{\mathrm{THH}(\mathbb{F}_p)} \mathrm{THH}(\mathbb{F}_p)[2])^{C_{p^r}}$ is a $\mathrm{THH}(\mathbb{F}_p)^{C_{p^r}}$-module for every $\mathrm{THH}(\mathbb{F}_p)$-module $N$ in pre-$p$-cyclotomic spectra, thus the graded abelian group

$$\pi_*(N \otimes^{\mathbb{L}}_{\mathrm{THH}(\mathbb{F}_p)} \mathrm{THH}(\mathbb{F}_p)[2])^{C_{p^r}}$$

is a $p^{r+1}$-torsion. Applying $\mathrm{TR}(M \otimes^{\mathbb{L}}_{\mathrm{THH}(\mathbb{F}_p)} (-))$ to the fiber sequence (3.1), we get a fiber sequence

$$0 = \mathrm{TR}(M \otimes^{\mathbb{L}}_{\mathrm{THH}(\mathbb{F}_p)} \mathrm{THH}(\mathbb{F}_p)[2]) \longrightarrow \mathrm{TR}(M) \longrightarrow \mathrm{TR}(M \otimes^{\mathbb{L}}_{\mathrm{THH}(\mathbb{F}_p)} \underline{W}(\mathbb{F}_p))$$

---

3.1. When $M = \mathrm{THH}(R)$, the restriction maps coming from the cyclonic structure are referred to as Frobenius maps, not the restriction maps coming from the pre-cyclotomic structure.

of spectra, and the result follows. □

*Remark* 3.3. In Proposition 3.2, if we only need to prove that the map in question becomes an equivalence after $p$-completion, we do not have to pass to homotopy groups. Indeed, after modulo $p$, the transition maps in (3.3) are not only phantom but zero, since the composite map (3.2) is zero after modulo $p$. As a consequence, the $p$-completion of $\mathrm{TR}(\mathrm{THH}(\mathbb{F}_p)[2])$ is zero as well, where $\mathrm{THH}(\mathbb{F}_p)[2]$ is equipped with the pre-$p$-cyclotomic structure as in the proof of Proposition 3.2. This argument works in slightly more general contexts, e.g. in non-Postnikov-complete topoi.

*Remark* 3.4. A similar argument, with slightly more efforts, proves a relative, mixed characteristic version of Proposition 3.2. We will address this in the future.

## 4. Relative Hochschild–Witt homology

In this section, we first define relative Hochschild–Witt homology, and then briefly explain that Hesselholt's proof of his HKR theorem in [Hes96] leads to an HKR theorem for relative Hochschild–Witt homology.

**Definition 4.1.** *Let $k$ be a commutative ring, and $\mathcal{C}$ a dualizable stable $k$-linear $\infty$-category. Then the* Hochschild–Witt homology $\underline{W}\,\mathrm{HH}(\mathcal{C}/k)$ *relative to $k$ is defined to be the topological Hochschild homology*

$$\mathrm{THH}(\mathcal{C}) \otimes^{\mathbb{L}}_{\mathrm{THH}(k)} \underline{W}(k)$$

*of $\mathcal{C}$ relative to the pre-$p$-cyclotomic base $\underline{W}(k)$. Let $r \in \mathbb{N}_{>0}$. The $r$-truncated Hochschild–Witt homology $W_r\,\mathrm{HH}(\mathcal{C}/k)$ relative to $k$ is defined to be the genuine $C_{p^{r-1}}$-fixed points*

$$(\underline{W}\,\mathrm{HH}(\mathcal{C}/k))^{C_{p^{r-1}}}$$

*of the relative Hochschild–Witt homology.*

It follows directly from Proposition 3.2 that

**Corollary 4.2.** *Let $\mathcal{C}$ be a dualizable stable $\mathbb{F}_p$-linear $\infty$-category. Then we have an equivalence*

$$\mathrm{TR}(\mathcal{C}) \longrightarrow \mathrm{TR}(\underline{W}\,\mathrm{HH}(\mathcal{C}/\mathbb{F}_p))$$

*of topological Cartier modules.*

In particular, combining with Hesselholt's HKR theorem, we see that, for a smooth commutative $\mathbb{F}_p$-algebra, we can identify Hochschild–Witt homology groups relative to $\mathbb{F}_p$ with de Rham–Witt forms. After identification with Kaledin's Hochschild–Witt homology, this recovers the limit case of Kaledin's HKR theorem. We now generalize this to the relative situation. We first construct an $F$-$V$-procomplex structure on $\pi_*(W_r\,\mathrm{HH}(A/k))_{r \in \mathbb{N}_{>0}}$ for commutative $k$-algebras $A$.

Let $k \to A$ be a map of commutative rings. Recall that, as in [LZ04, Def 1] (reformulated in [BMS18, Def 10.5]), an *$F$-$V$-procomplex* over $k \to A$ consists of the following data:

1. for every $r \in \mathbb{N}_{>0}$, a commutative differential ($\mathbb{N}$-)graded $W_r(k)$-algebra (abbrev. $W_r(k)$-CDGA) $P_r^{\cdot} = (P_r^n)_{n \in \mathbb{N}}$.

2. for every $r \in \mathbb{N}_{>0}$, a map $F \colon P_{r+1}^{\cdot} \to F_* P_r^{\cdot}$ of graded $W_{r+1}(k)$-algebras (not CDGA's!), where the second $F$ is the Frobenius $F \colon W_{r+1}(k) \to W_r(k)$;

3. for every $r \in \mathbb{N}_{>0}$, a map $V: F_* P_r^\cdot \to P_{r+1}^\cdot$ of graded $W_{r+1}(k)$-modules;
4. for every $r \in \mathbb{N}_{>0}$, a map $\lambda_r: W_r(A) \to P_r^0$ of $W_r(k)$-algebras, compatible with the maps $F$, $V$, and $R$,
5. for every $r \in \mathbb{N}_{>0}$, a map $R: P_{r+1}^\cdot \to R_* P_r^\cdot$ of $W_{r+1}(k)$-CDGA's, where the second $R$ is the restriction map $R: W_{r+1}(k) \to W_r(k)$;

which satisfies the following properties:

1. $FV = p$;
2. $F \, \mathrm{d} V = \mathrm{d}$;
3. Frobenius reciprocity: $V(F(x) \, y) = x \, V(y)$;
4. Teichmüller identity: $F \, \mathrm{d} \lambda_{r+1}([x]) = \lambda_r([x])^{p-1} \, \mathrm{d} \lambda_r([x])$ for $x \in A$ and $r \in \mathbb{N}_{>0}$.
5. $RF = FR$, $RV = VR$;

*Construction* 4.3. **(Hesselholt)** Let $k \to A$ be a map of commutative rings. Then the sequence $(\pi_\cdot W_r \operatorname{HH}(A/k))_{r \in \mathbb{N}_{>0}}$ carries an $F$-$V$-procomplex (over $k \to A$) structure.

Indeed, given any $\mathbb{E}_\infty$-$\underline{W}(k)$-algebra $B$ in $\operatorname{Sp}^{g_p B \mathbb{T}}$, the homotopy groups $\left( \pi_\cdot B^{C_{p^{r-1}}} \right)_{r \in \mathbb{N}_{>0}}$ has the first three items of data:

1. for each $r \in \mathbb{N}_{>0}$, the $\mathbb{T}/C_{p^{r-1}}$-equivariant $\mathbb{E}_\infty$-$W_r(k)$-structure on $B^{C_{p^{r-1}}}$ gives rise to a $W_r(k)$-CDGA structure on $\left( \pi_n B^{C_{p^{r-1}}} \right)_{n \in \mathbb{N}}$;
2. for each $r \in \mathbb{N}_{>0}$, the map $F: (\pi_n B^{C_{p^r}}) \to \left( F_* \pi_n B^{C_{p^{r-1}}} \right)$ is induced by the restriction $\mathbb{E}_\infty$-map $B^{C_{p^r}} \to B^{C_{p^{r-1}}}$, which is over the Frobenius map $W_{r+1}(k) = \underline{W}(k)^{C_{p^r}} \to \underline{W}(k)^{C_{p^{r-1}}} = W_r(k)$;
3. for each $r \in \mathbb{N}_{>0}$, the map $V: \left( F_* \pi_n B^{C_{p^{r-1}}} \right) \to (\pi_n B^{C_{p^r}})$ is induced by the transfer map $B^{C_{p^{r-1}}} \to B^{C_{p^r}}$, and the first three properties are satisfied: the first and the last are direct consequences of graded (anti-commutative) Green functors; the second follows from [Hes96, Lem 1.5.1] (all cyclonic spectra here are over $\underline{W}(\mathbb{Z})$, and there the Hopf element $\eta$ vanishes).

The $\lambda$-map and the Teichmüller identity essentially follows from the norm structure on $\underline{W} \operatorname{HH}(A/k)$. Here we give an elementary argument: we know that

$$\begin{aligned} \pi_0 \underline{W} \operatorname{HH}(A/k) &= \pi_0(\operatorname{THH}(A) \otimes^{\mathbb{L}}_{\operatorname{THH}(k)} \pi_0 \operatorname{THH}(k)) \\ &= \pi_0(\operatorname{THH}(A)) \\ &= \underline{W}(A) \end{aligned}$$

which gives rise to the maps $\lambda_r$ (which are even isomorphisms). Moreover, the Teichmüller identity follows directly from [Hes96, Lem I.5.6], and the fact that the canonical map $(\pi_\cdot \operatorname{TR}^r(A))_{r \in \mathbb{N}_{>0}} \to (\pi_\cdot W_r \operatorname{HH}(A/k))_{r \in \mathbb{N}_{>0}}$ preserves the maps $F$, $d$, $\lambda$.

Finally, the map $R$ follows from the map $R: W_{r+1} \operatorname{HH}(A/k) \to W_r \operatorname{HH}(A/k)$ as in Construction 2.6, and the compatibility follows from the $\mathbb{E}_\infty$-cyclonic ring structure on $R$.

*Construction* 4.4. Let $k \to A$ be a map of commutative rings. Then the universal property of relative de Rham–Witt complexes in [LZ04] induces a unique map

$$(W_r \Omega^{\cdot}_{A/k})_{r \in \mathbb{N}_{>0}} \longrightarrow (\pi_* W_r \operatorname{HH}(A/k))_{r \in \mathbb{N}_{>0}}$$

of $F$-$V$-procomplexes.

Now the proof in [Hes96, §2] (where the key is the computational results in Prop 2.2.5 & 2.3.3) implies that

**Lemma 4.5.** *Let $k$ be a commutative ring, and $A$ a finite polynomial $k$-algebra. Then the map*

$$(W_r \Omega^{\cdot}_{A/k})_{r \in \mathbb{N}_{>0}} \longrightarrow (\pi_* W_r \operatorname{HH}(A/k))_{r \in \mathbb{N}_{>0}}$$

*in Construction 4.4 is an isomorphism.*

*Remark* 4.6. Let $A \to B$ be an étale map of commutative $k$-algebras, and $r \in \mathbb{N}_{>0}$. Then by [Bor11b, 15.4], the map

$$\underline{W}(A) \otimes^{\mathbb{L}}_{W_r(A)} W_r(B) \longrightarrow \underline{W}(B)$$

in $D(W_r(A)) \otimes \operatorname{Sp}^{gBC_{p^{r-1}}}$ induced by the map $\underline{W}(A) \to \underline{W}(B)$ in $D(W_r(B)) \otimes \operatorname{Sp}^{gBC_{p^{r-1}}}$ is an equivalence. Furthermore, by [Bor11a, Thm B], the map $W_r(A) \to W_r(B)$ is étale. It follows that the map $\underline{W}(A) \to \underline{W}(B)$ of $\mathbb{E}_\infty$-cyclonic rings is étale as well.

It follows from Remark 4.6 and [HLL20, Add 3.2] that

**Lemma 4.7.** *Let $A \to B$ be an étale map of commutative $k$-algebras. Then the map*

$$\operatorname{THH}(A/k) \longrightarrow \operatorname{THH}(B/k)$$

*in $\operatorname{CAlg}(\operatorname{Sp}^{g_p B\mathbb{T}})$ is flat (and even étale).*

**Corollary 4.8.** *Let $k$ be a commutative ring, and $A \to B$ an étale map of commutative $k$-algebras. Then the map*

$$\underline{W} \operatorname{HH}(A/k) \longrightarrow \underline{W} \operatorname{HH}(B/k)$$

*in $\operatorname{CAlg}(\operatorname{Sp}^{g_p B\mathbb{T}})$ is flat (and even étale).*

We are ready to establish the HKR theorem:

**Proposition 4.9.** *Let $k$ be a commutative ring, and $A$ a smooth commutative $k$-algebra. Then the map*

$$(W_r \Omega^{\cdot}_{A/k})_{r \in \mathbb{N}_{>0}} \longrightarrow (\pi_* W_r \operatorname{HH}(A/k))_{r \in \mathbb{N}_{>0}}$$

*in Construction 4.4 is an isomorphism.*

*Proof.* By [BMS18, Lem 10.8], given an étale map $B \to C$ of $k$-algebras, the canonical map

$$W_r \Omega^{\cdot}_{B/k} \otimes_{W_r(B)} W_r(C) \longrightarrow W_r \Omega^{\cdot}_{C/k}$$

is an isomorphism, and by Corollary 4.8, the canonical map

$$\pi_* W_r \operatorname{HH}(B/k) \otimes_{W_r(B)} W_r(C) \longrightarrow \pi_* W_r \operatorname{HH}(C/k)$$

is an isomorphism as well. It follows that, the map in Construction 4.4 is an isomorphism for étale extensions of polynomial $k$-algebras. Now by smoothness of $A$, there exists a (finite) Zariski cover $A \to \prod_{i \in I} A_i$ such that each $A_i$ is an étale extension of a polynomial $k$-algebra, where $I$ is a finite set. Then we examine the commutative diagram

$$\begin{array}{ccc} W_r\Omega^{\cdot}_{A/k} \otimes_{W_r(A)} \left( \prod_{i \in I} W_r(A_i) \right) & \longrightarrow & \pi_* W_r \operatorname{HH}(A/k) \otimes_{W_r(A)} \left( \prod_{i \in I} W_r(A_i) \right) \\ \downarrow & & \downarrow \\ \prod_{i \in I} W_r\Omega^{\cdot}_{A_i/k} & \longrightarrow & \prod_{i \in I} \pi_* W_r \operatorname{HH}(A_i/k) \end{array}$$

where the bottom horizontal arrow and the vertical arrows are isomorphisms, thus the top horizontal arrow is an isomorphism as well. Now the result follows from [Bor11a, Prop 6.9], which implies that the map $W_r(A) \to \prod_{i \in I} W_r(A_i)$ is faithfully flat. □

## 5. The norm over $\underline{\mathbb{Z}}$

Let $G$ be a finite group, and $X \in D(\mathbb{Z})^{BG}$ a $G$-equivariant object in $D(\mathbb{Z})$. Then we can talk about *homotopy orbits* and *homotopy fixed points* of $X$. However, sometimes one needs a stricter notion of orbits and fixed points. For example, let $d \in \mathbb{N}$ be a natural number. Then we have the $d$-th derived symmetric power

$$\operatorname{LSym}^d_{\mathbb{Z}} \colon D(\mathbb{Z}) \longrightarrow D(\mathbb{Z})$$

obtained by right-left extending ([BM19, §3.2], see also [Rak20, §4.2]) the "usual" symmetric power functor

$$\begin{array}{rcl} \operatorname{Latt}_{\mathbb{Z}} & \longrightarrow & D(\mathbb{Z}) \\ M & \longmapsto & (M^{\otimes d})_{\Sigma_d} \end{array}$$

polynomial of degree $d$, defined on the category $\operatorname{Latt}_{\mathbb{Z}}$ of finite free abelian groups. It would be desirable if, for every $\mathbb{Z}$-module spectrum $M \in D(\mathbb{Z})$, we can define $M^{\otimes d}$ in a world where we can talk about stricter $\Sigma_d$-orbits. In [BCN21, Ex 2.15 & 2.51 & 3.68], the authors showed that the ∞-category $\operatorname{Mod}_{\underline{\mathbb{Z}}}(\operatorname{Sp}^{\Sigma_d})$ of *derived cohomological $\Sigma_d$-Mackey functors* is such a world: we have functors

$$\begin{array}{rcl} (-)^{\otimes d} \colon D(\mathbb{Z}) & \longrightarrow & \operatorname{Mod}_{\underline{\mathbb{Z}}}(\operatorname{Sp}^{\Sigma_d}) \\ (-)_{\Sigma_d} \colon \operatorname{Mod}_{\underline{\mathbb{Z}}}(\operatorname{Sp}^{\Sigma_d}) & \longrightarrow & D(\mathbb{Z}) \end{array}$$

where the first functor is a sifted-colimit-preserving polynomial functor of degree $d$, and the second functor preserves small colimits. In this section, we will show that the first functor $(-)^{\otimes d}$ is in fact a special case of norms relative[5.1] to $\underline{\mathbb{Z}}$. For this, we first recall that, letting $X$ be a finite groupoid, then a *finite $X$-set* [BH21, §9.1] is a functor $X \to \operatorname{Fin}$, and an *$X$-module* is a functor $X \to \operatorname{Ab}$.

**Definition 5.1.** *Let $X$ be a finite groupoid*[5.2]. *A* finite permutation $X$-module *is an $X$-module freely generated by a finite $X$-set, and the category of finite permutation $X$-modules is denoted by* $\operatorname{Perm}_X$.

---

[5.1]. We learnt this notion from Kaif Hilman, and this observation was formulated during a discussion with him.

[5.2]. A groupoid is *finite* if it is finite as a 1-category (not as an ∞-category in general).

In [BCN21, Ex 2.15], the authors showed that we have an equivalence
$$\mathrm{Mod}_{\underline{\mathbb{Z}}}(\mathrm{Sp}^{gX}) \simeq \mathrm{Sp}(\mathcal{P}_\Sigma(\mathrm{Perm}_X)) =: D\,\mathrm{Mack}^{\mathrm{coh}}_X$$
of $\infty$-categories when $X$ is a groupoid. We now upgrade[5.3] it to an equivalence between the global symmetric monoidal $\infty$-category $\mathrm{Mod}_{\underline{\mathbb{Z}}}(\mathrm{Sp}^{g(-)}): X \mapsto \mathrm{Mod}_{\underline{\mathbb{Z}}}(\mathrm{Sp}^{gX})$ of $\underline{\mathbb{Z}}$-modules and the global symmetric monoidal $\infty$-category $D\,\mathrm{Mack}^{\mathrm{coh}}_{(-)}$ of derived cohomological Mackey functors, with the following definition of global symmetric monoidal categories.

**Definition 5.2.** *A global symmetric monoidal $\infty$-category $\mathcal{C}$ is a finite-product-preserving functor*
$$\mathcal{C}^\otimes : \mathrm{Span}(\mathrm{Gpd}^{\mathrm{fin}}, \mathrm{all}, \mathrm{faith}) \longrightarrow \mathrm{Cat}_\infty.$$
*We will denote $\mathcal{C}^\otimes(BG)$ by $\mathcal{C}_G$ for every finite group $G$.*

We start with the construction of the global symmetric monoidal categories $\mathrm{Perm}_{(-)}$ and $D\,\mathrm{Mack}^{\mathrm{coh}}_{(-)}$.

*Remark* 5.3. Let $\mathcal{C}$ be an $\infty$-category. It is known that there is a canonical global symmetric monoidal $\infty$-category
$$\begin{aligned} \mathcal{C}^{(-)} : \mathrm{Span}(\mathrm{Gpd}^{\mathrm{fin}}, \mathrm{all}, \mathrm{faith}) &\longrightarrow \mathrm{Cat}_1 \\ X &\longmapsto \mathcal{C}^X. \end{aligned}$$
Apply this to $\mathcal{C} = \mathrm{Ab}$, we obtain the global symmetric monoidal (1-)category $\mathrm{Ab}^{(-)}$, which is globally presentable and additive.

*Construction* 5.4. The global symmetric monoidal (1-)category $\mathrm{Ab}^{(-)}$ admits a full global subcategory $\mathrm{Perm}_{(-)} \subseteq \mathrm{Ab}^{(-)}$ which is additive and inherits a globally additive symmetric monoidal structure as well. Consequently, we get another globally presentably stable symmetric monoidal $\infty$-category $D\,\mathrm{Mack}^{\mathrm{coh}}_{(-)} := \mathrm{Sp}(\mathcal{P}_\Sigma(\mathrm{Perm}^{(-)}))$ obtained by taking fiberwise stablized non-abelian derived categories.

*Construction* 5.5. Note that the unit in $D\,\mathrm{Mack}^{\mathrm{coh}}_{(-)}$ is globally compact by construction, thus we get a globally presentable adjunction
$$\mathrm{Sp}^{g(-)} \rightleftarrows D\,\mathrm{Mack}^{\mathrm{coh}}_{(-)}. \tag{5.1}$$
In particular, the unit in $D\,\mathrm{Mack}^{\mathrm{coh}}_{(-)}$ defines a global $\mathbb{E}_\infty$-ring, denoted by $\underline{\mathbb{Z}}$ (we will soon justify this notation), via the global lax symmetric monoidal right adjoint $D\,\mathrm{Mack}^{\mathrm{coh}}_{(-)} \to \mathrm{Sp}^{g(-)}$, and we obtain a globally presentable adjunction
$$\mathrm{Mod}_{\underline{\mathbb{Z}}}(\mathrm{Sp}^{g(-)}) \rightleftarrows D\,\mathrm{Mack}^{\mathrm{coh}}_{(-)}. \tag{5.2}$$

It follows from the equivariant symmetric monoidal Schwede–Shipley theorem that

**Proposition 5.6.** *The adjunction (5.2) is an equivalence of globally presentably stable symmetric monoidal $\infty$-categories.*

Now we show that the global $\mathbb{E}_\infty$-ring $\underline{\mathbb{Z}}$ coincides with the constant global Tambara functor $\underline{\mathbb{Z}}$. We first show that it lies in the heart.

---

[5.3]. This formulation and basic ideas of the argument were proposed by Marc HOYOIS. We also thank Kaif HILMAN and Sil LINSKENS for discussions.

**Lemma 5.7.** *Let $X$ be a finite groupoid. Then the spectrum $\underline{\mathbb{Z}}(X) \in \mathrm{Sp}^{gX}$ lies in the heart $\mathrm{Sp}^{gX,\heartsuit} = \mathrm{Mack}_X(\mathrm{Ab})$. The global $\mathbb{E}_\infty$-ring $\underline{\mathbb{Z}}$ is (the Eilenberg–MacLane spectrum of) a global Tambara functor, whose underlying Mackey functor is the constant Mackey functor $\underline{\mathbb{Z}}$.*

*Proof.* For every finite $X$-set $S$, by adjunction (5.1), the anima

$$\begin{aligned}\mathrm{Hom}_{\mathrm{Sp}^{gX}}(\Sigma_X^\infty S_+, \underline{\mathbb{Z}}(X)) &= \mathrm{Hom}_{D\mathrm{Mack}_X^{\mathrm{coh}}}(\mathbb{Z}[S], \mathbb{Z}) \\ &= \mathrm{Hom}_{\mathrm{Perm}_X}(\mathbb{Z}[S], \mathbb{Z})\end{aligned}$$

is a set. □

*Remark* 5.8. The adjunction (5.1) induces a globally presentably additive adjunction

$$\mathrm{Mack}_{(-)}(\mathrm{Ab}) \rightleftarrows \mathrm{Fun}^\oplus(\mathrm{Perm}_X^{\mathrm{op}}, \mathrm{Ab}) =: \mathrm{Mack}_{(-)}^{\mathrm{coh}} \tag{5.3}$$

on the hearts.

**Corollary 5.9.** *The global $\mathbb{E}_\infty$-ring $\underline{\mathbb{Z}}$ can be identified with the image of the unit in $\mathrm{Mack}_{(-)}^{\mathrm{coh}}$ under the right adjoint $\mathrm{Mack}_{(-)}^{\mathrm{coh}} \to \mathrm{Mack}_{(-)}(\mathrm{Ab})$.*

*Remark* 5.10. The canonical inclusion $\mathrm{Perm}_{(-)} \subseteq \mathrm{Ab}^{(-)}$ induces an globally additive adjunction

$$\mathrm{Mack}_{(-)}^{\mathrm{coh}} \rightleftarrows \mathrm{Ab}^{(-)}.$$

Moreover, the right adjoint carries the unit in $\mathrm{Ab}^{(-)}$ to the unit in $\mathrm{Mack}_{(-)}^{\mathrm{coh}}$, since it is a unit in $\mathrm{Perm}_{(-)}$. Consequently, we can identify the global $\mathbb{E}_\infty$-ring $\underline{\mathbb{Z}}$ with the image of the unit in $\mathrm{Ab}^{(-)}$ under the right adjoint $\mathrm{Ab}^{(-)} \to \mathrm{Mack}_{(-)}(\mathrm{Ab})$.

Now the identification of the global $\mathbb{E}_\infty$-ring $\underline{\mathbb{Z}}$ with the constant global Tambara functor $\underline{\mathbb{Z}}$ follows from the following general statement.

**Lemma 5.11.** *Let $k$ be a commutative ring, viewed as a global $\mathbb{E}_\infty$-algebra in the global symmetric monoidal category $\mathrm{Ab}^{(-)}$. Then its image under the right adjoint $G \colon \mathrm{Ab}^{(-)} \to \mathrm{Mack}_{(-)}(\mathrm{Ab})$ can be identified with the constant global Tambara functor $\underline{k}$.*

*Proof.* We compute the norm along an inclusion $H \subseteq H'$ of groups as an example, and this is a special case of [Hil22, Thm 2.4.10]: the norm on $G(k)$ is given by the map

$$\pi_0(N_H^{H'} \underline{k}) \longrightarrow \underline{k}^{\otimes(H'/H)} \longrightarrow \underline{k}$$

of $H$-Mackey functors, which corresponds to the $|H'/H|$-th power map $G(k)(BH)^H \to G(k)(BH)^{H'}$, i.e. the norm map $\underline{k}(BH)^H \to \underline{k}(BH')^{H'}$. □

*Example* 5.12. Let $n \in \mathbb{N}_{>0}$ be a positive integer. Then the norm $i_\otimes$ of the global symmetric monoidal category $\mathrm{Ab}^{(-)}$ along the inclusion $i \colon \{e\} \hookrightarrow C_n$ of finite groups is concretely given by the tensor power functor

$$\begin{aligned}\mathrm{Ab} &\longrightarrow \mathrm{Ab}^{BC_n} \\ M &\longmapsto M^{\otimes n},\end{aligned}$$

where $C_n$ acts on $M^{\otimes n}$ by permuting factors, and consequently, the norm $i_\otimes$ of the global symmetric monoidal $\infty$-category $D\mathrm{Mack}_{(-)}$ is the right-left-extension of this tensor power functor on finite free abelian groups. We will denote the norm $i_\otimes$ by $(-)^{\otimes_{\mathbb{Z}}^{\mathbb{L}} C_n}$ (or $(-)^{\otimes^{\mathbb{L}} C_n}$ when there is no ambiguity).

*Example* 5.13. More generally, let $G$ be a finite group, and $S$ a finite $G$-set. Then the norm of the global symmetric monoidal category $\mathrm{Ab}^{(-)}$ along the span $* \leftarrow S_{hG} \to BG$[5.4] of finite groupoids is given by

$$\begin{aligned} \mathrm{Ab} &\longrightarrow \mathrm{Ab}^G \\ M &\longmapsto M^{\otimes S}, \end{aligned}$$

where $G$ acts on $M^{\otimes S}$ by permuting factors, and consequently, the norm of the global symmetric monoidal $\infty$-category $D\,\mathrm{Mack}_{(-)}$ along the span $* \leftarrow S_{hG} \to BG$ is the right-left-extension of this tensor power functor on finite free abelian groups. In particular, when $G = \Sigma_n$ and $S = \{1, \ldots, n\}$ where $G$ acts canonically, we recover a simplified version of [BCN21, Rem 3.68].

## 6. The norm over $\underline{W}(\mathbb{F}_p)$

In this section, we relate Kaledin's *polynomial Witt vectors* [Kal18] to norms over the ($p$-typical) $\mathbb{T}$-$\mathbb{E}_\infty$-$\underline{\mathbb{Z}}$-algebra $\underline{W}(\mathbb{F}_p)$, whose definition will be reviewed below. More precisely, we identify Kaledin's Mackey-valued polynomial Witt vectors $\tilde{W}_r(M)$ with the norm $M^{\otimes^{\mathbb{L}}_{\underline{W}(\mathbb{F}_p)} C_{p^{r-1}}}$ relative to $\underline{W}(\mathbb{F}_p)$, functorially in $\mathbb{F}_p$-vector spaces $M$ (Proposition 6.5). Consequently, Kaledin's polynomial Witt vector $W_r(M)$ admits an "explicit" formula $\left(M^{\otimes^{\mathbb{L}}_{\underline{W}(\mathbb{F}_p)} C_{p^{r-1}}}\right)^{C_{p^{r-1}}}$ (Corollary 6.6). In fact, we establish a parametrized version of this formula (Lemma 6.4), which allows us to identify the polynomial Witt trace theory in Section 7.

We start with an alternative construction of Kaledin's Mackey-valued Tate-construction $Q$.

*Construction* 6.1. Let $G$ be a finite group. The constant map $G \to *$ of finite $G$-sets induces the augmentation map $\mathrm{aug}_G: \mathbb{Z}[G] \to \mathbb{Z}$ of permutation $G$-modules, and its cokernel $\mathrm{coker}(\mathrm{aug}_G)$, taken in the abelian category $\mathrm{Mack}^{\mathrm{coh}}_{BG}$, gives rise to an endofunctor

$$(-) \otimes^{\mathbb{L}}_{\underline{\mathbb{Z}}} \mathrm{coker}(\underline{\mathrm{aug}}_G) : D\,\mathrm{Mack}^{\mathrm{coh}}_{BG} \longrightarrow D\,\mathrm{Mack}^{\mathrm{coh}}_{BG}.$$

**Lemma 6.2.** *Let $G$ be a finite group. Then the endofunctor $(-) \otimes^{\mathbb{L}}_{\underline{\mathbb{Z}}} \mathrm{coker}(\underline{\mathrm{aug}}_G)$ in Construction 6.1 coincides with Kaledin's functor $Q$ in [Kal18, Def 1.3] after restriction to permutation $G$-modules.*

*Proof.* Since permutation $G$-modules are flat in $\mathrm{Mack}^{\mathrm{coh}}_{BG}$, it suffices to identify the non-derived composite functor

$$(-) \otimes_{\underline{\mathbb{Z}}} \mathrm{coker}(\underline{\mathrm{aug}}_G) : \mathrm{Perm}_{BG} \hookrightarrow \mathrm{Mack}^{\mathrm{coh}}_{BG} \longrightarrow \mathrm{Mack}^{\mathrm{coh}}_{BG}$$

with Kaledin's functor $Q$ in [Kal18, Def 1.3]. By definition, the functor $(-) \otimes_{\underline{\mathbb{Z}}} \underline{\mathbb{Z}[G]}$ coincides with Kaledin's $\psi := \psi^{\{e\}}$, and the natural transformation

$$\underline{\mathrm{aug}_G} : (-) \otimes_{\underline{\mathbb{Z}}} \underline{\mathbb{Z}[G]} \longrightarrow (-)$$

induced by the augmentation map $\mathrm{aug}_G : \mathbb{Z}[G] \to \mathbb{Z}$ factors through $\psi(-)_G$, thus $\mathrm{coker}(\underline{\mathrm{aug}_G}) = \mathrm{coker}(\mathrm{tr})$ where $\mathrm{tr}: \psi(-)_G \to (-)$ is the natural transformation [Kal18, (1.25)]. □

---

5.4. More generally, given two finite groups $H$, $G$, and an $H$-$G$-biset $S$, we can form a span $BH \leftarrow [H \backslash S / G] \to BG$ of finite groupoids. We learned this generality from Marc Hoyois.

When $G$ is a finite cyclic $p$-group, the cokernel $\operatorname{coker}(\operatorname{aug}_G)$ can be related to the cohomological $\mathbb{T}$-Mackey functor $\underline{W}(\mathbb{F}_p)$. More precisely,

**Lemma 6.3.** *Let $r \in \mathbb{N}_{>0}$ be a positive integer. Then we have identifications*
$$\underline{W}(\mathbb{F}_p) \simeq \operatorname{coker}(p \operatorname{aug}_{C_{p^r}})$$
$$\operatorname{Infl}_{C_{p^{r-1}}}^{C_{p^r}} \underline{W}(\mathbb{F}_p) \simeq \operatorname{coker}(\operatorname{aug}_{C_{p^r}})$$
*in $D\operatorname{Mack}_{C_{p^{r-1}}}^{\operatorname{coh}}$. In other words, Kaledin's functor $Q$ on finite cyclic p-groups can be read off from the $\mathbb{T}$-symmetric monoidal functor*
$$\operatorname{Mod}_{\underline{\mathbb{Z}}}(\operatorname{Sp}^{g_p B\mathbb{T}}) \xrightarrow{(-) \otimes_{\underline{\mathbb{Z}}}^{\mathbb{L}} \underline{W}(\mathbb{F}_p)} \operatorname{Mod}_{\underline{W}(\mathbb{F}_p)}(\operatorname{Sp}^{g_p B\mathbb{T}}).$$

*Proof.* It follows from an explicit computation. In fact, [Zen17, Ex 3.19] gives an explicit resolution
$$0 \longrightarrow \underline{\mathbb{Z}} \longrightarrow \underline{\mathbb{Z}}[C_{p^{r-1}}] \longrightarrow \underline{\mathbb{Z}}[C_{p^{r-1}}] \xrightarrow{p \operatorname{aug}_{C_{p^{r-1}}}} \underline{\mathbb{Z}} \longrightarrow \underline{W}(\mathbb{F}_p) \longrightarrow 0$$
of $\underline{W}(\mathbb{F}_p)$ in the abelian category $\operatorname{Mack}_{C_{p^{r-1}}}^{\operatorname{coh}}$. □

The comparison with Kaledin's polynomial Witt vectors essentially follows from combining Lemma 6.3 and the following Lemma 6.4, which we will explain soon.

**Lemma 6.4.** *Let $f: X \to Y$ be a map in $\operatorname{Fin}_{B\mathbb{T}}$, viewed as a map of finite groupoids over $B\mathbb{T}$. Then we have a commutative diagram*
$$\begin{array}{ccc} \operatorname{Mod}_{\underline{\mathbb{Z}}}(\operatorname{Sp}^{gX}) & \xrightarrow{f_{\otimes}} & \operatorname{Mod}_{\underline{\mathbb{Z}}}(\operatorname{Sp}^{gY}) \\ \downarrow & & \downarrow \\ \operatorname{Mod}_{\underline{W}(\mathbb{F}_p)}(\operatorname{Sp}^{gX}) & \xrightarrow{f_{\otimes}} & \operatorname{Mod}_{\underline{W}(\mathbb{F}_p)}(\operatorname{Sp}^{gY}) \end{array}.$$
*In particular, letting $r \in \mathbb{N}_{>0}$ be a positive integer, and taking $X = *$ and $Y = BC_{p^{r-1}}$, we get a commutative diagram*
$$\begin{array}{ccc} D(\underline{\mathbb{Z}}) & \xrightarrow{(-)^{\otimes_{\underline{\mathbb{Z}}}^{\mathbb{L}} C_{p^{r-1}}}} & D\operatorname{Mack}_{BC_{p^{r-1}}}^{\operatorname{coh}} \\ \downarrow & & \downarrow \\ D(\mathbb{F}_p) & \xrightarrow{(-)^{\otimes_{\underline{W}(\mathbb{F}_p)}^{\mathbb{L}} C_{p^{r-1}}}} & \operatorname{Mod}_{\underline{W}(\mathbb{F}_p)}(\operatorname{Sp}^{gBC_{p^{r-1}}}) \end{array}.$$

Now we restrict to the full subcategories $\operatorname{Mod}_{\mathbb{Z}}^{\operatorname{free}} \subseteq D(\mathbb{Z})$ and $\operatorname{Mod}_{\mathbb{F}_p}^{\operatorname{free}} \subseteq D(\mathbb{F}_p)$ (not necessarily finite, but functors in question do preserve filtered colimits). Note that the base change functor $\operatorname{Mod}_{\mathbb{Z}}^{\operatorname{free}} \to \operatorname{Mod}_{\mathbb{F}_p}^{\operatorname{free}}$ is essentially surjective, and epimorphic on Hom's. Combining with [Kal18, Prop 2.3], we get

**Proposition 6.5.** *Let $r \in \mathbb{N}_{>0}$ be a positive integer. Then there is a commutative diagram*
$$\begin{array}{ccc} \operatorname{Mod}_{\mathbb{F}_p}^{\operatorname{free}} & \xrightarrow{\tilde{W}_r} \operatorname{Mack}_{BC_{p^{r-1}}}(\operatorname{Ab}) & \hookrightarrow \operatorname{Sp}^{gBC_{p^{r-1}}} \\ \downarrow & & \uparrow \\ D(\mathbb{F}_p) & \xrightarrow{(-)^{\otimes_{\underline{W}(\mathbb{F}_p)}^{\mathbb{L}} C_{p^{r-1}}}} & \operatorname{Mod}_{\underline{W}(\mathbb{F}_p)}(\operatorname{Sp}^{gBC_{p^{r-1}}}) \end{array}$$

of $\infty$-categories. In particular, Mackey-valued polynomial Witt vectors are cohomological Mackey functors.

**Corollary 6.6.** *Let $r \in \mathbb{N}_{>0}$. The composite functor*

$$\mathrm{Mod}_{\mathbb{F}_p}^{\mathrm{free}} \hookrightarrow D(\mathbb{F}_p) \xrightarrow{(-)^{\otimes_{\underline{W}(\mathbb{F}_p)}^{\mathbb{L}} C_{p^{r-1}}}} \mathrm{Mod}_{\underline{W}(\mathbb{F}_p)}(\mathrm{Sp}^{gBC_{p^{r-1}}}) \xrightarrow{(-)^{C_{p^{r-1}}}} D(W_r(\mathbb{F}_p))$$

*coincides with $r$-truncated polynomial Witt vectors.*

## 7. From norms to trace theories

Let $\mathcal{C}$ be a symmetric monoidal[7.1] 1-category, and $\mathcal{E}$ a 1-category. Recall that an $\mathcal{E}$-valued trace theory on $\mathcal{C}$, as defined in [Kal15, Def 2.1][7.2], is a functor $T: \mathcal{C} \to \mathcal{E}$ equipped with a natural isomorphism $(\tau_{M,N})_{(M,N) \in \mathcal{C}^2}$ between functors $(M,N) \mapsto F(M \otimes N)$ and $(M,N) \mapsto F(N \otimes M)$, such that

**Unity.** for every $M \in \mathcal{C}$, we have $\tau_{1,M} = \mathrm{id}$, where we use unity equivalences $1 \otimes M \simeq M \simeq M \otimes 1$;

**Acyclicity.** for every triple $(L, M, N) \in \mathcal{C}^3$, we have

$$\tau_{L,M,N} \circ \tau_{N,L,M} \circ \tau_{M,N,L} = \mathrm{id},$$

where for $(A, B, C) \in \mathcal{C}^3$, we denote $\tau_{A,B,C} := \tau_{A, B \otimes C}$, and we use associativity equivalences $(A \otimes B) \otimes C \simeq A \otimes (B \otimes C)$.

The symmetric monoidal structure on $\mathcal{C}$ gives rise to a canonical $\mathcal{C}$-valued trace theory on $\mathcal{C}$, [Kal15, Ex 2.2]. A nontrivial important trace theory is given in [Kal15, Ex 2.4], which we recall as follows. Let $k$ be a commutative ring, $m \in \mathbb{N}_{>0}$ a positive integer, and we consider the category $\mathcal{C} = \mathrm{Proj}_k^{\mathrm{fg}}$ of finite projective $k$-modules. Then we have a $\mathrm{Mod}_k$-valued trace theory on $\mathcal{C}$ given by the functor $T: \mathcal{C} \to \mathrm{Mod}_k, M \mapsto (M^{\otimes s})_{C_m}$, the $C_m$-orbits of the $s$-th tensor power, with the natural isomorphism $\tau$ concretely induced by the isomorphism

$$\begin{aligned}
(M \otimes N)^{\otimes m} &= M \otimes N \otimes \cdots \otimes M \otimes N \\
&= M \otimes (N \otimes \cdots \otimes M \otimes N) \\
&\xrightarrow{\simeq} (N \otimes \cdots \otimes M \otimes N) \otimes M \\
&= N \otimes M \otimes \cdots \otimes N \otimes M \\
&= (N \otimes M)^{\otimes m}.
\end{aligned}$$

In fact, the functor $T$ factors as $\mathcal{C} \xrightarrow{(-)^{\otimes m}} \mathrm{Fun}(BC_m, \mathrm{Mod}_k) \xrightarrow{(-)} \mathrm{Mod}_k$. However, the first tensor power functor does not have a trace theory structure in general. Indeed, for trace theories on $\mathcal{C}$, one can verify that, for every pair $(M, N) \in \mathcal{C}^2$, we have $\tau_{M,N} \circ \tau_{N,M} = \mathrm{id}$ (as a consequence of unity and acyclicity), which is false in general unless $m = 1$. However, the tensor power admits a *subdivided trace theory*, a concept which is already implicit in [Kal15, Ex 2.7] and [Kal18] as well. In this section, we first review an $\infty$-categorical generalization of these concepts. Then we show that, as a generalization of the tensor power functor admitting a subdivided trace theory structure, the norm (or $C_m$-symmetric monoidal) structure gives rise to subdivided trace theories. We refer to [KMN23, §6.2] and [Nik] for relevant $\infty$-categorical accounts.

---

7.1. In fact, we can talk about trace theories on ($\mathbb{E}_1$-)monoidal $\infty$-categories. We do not need this generality.

7.2. Kaledin call them *trace functors*.

**Definition 7.1.** *Let $\mathcal{C}$ be a symmetric monoidal $\infty$-category, $m \in \mathbb{N}_{>0}$ a positive integer, and $\mathcal{E}$ an $\infty$-category. An $\mathcal{E}$-valued $m$-subdivided trace theory on $\mathcal{C}$ is a functor*
$$\Lambda_m^{\mathrm{op}} \times_{\mathrm{Fin}_*} \mathcal{C}^{\otimes} \longrightarrow \mathcal{E}$$
*which maps every coCartesian edges of the coCartesian fibration $\Lambda_m^{\mathrm{op}} \times_{\mathrm{Fin}_*} \mathcal{C}^{\otimes} \to \Lambda_m^{\mathrm{op}}$ to an equivalence in $\mathcal{E}$, where the functor $\Lambda_m^{\mathrm{op}} \to \mathrm{Fin}_*$ is the composite*
$$\Lambda_m^{\mathrm{op}} \xrightarrow{\mathrm{Cut}} \Lambda_m \xrightarrow{S \mapsto [S/\mathbb{Z}]} \mathrm{Fin} \xrightarrow{(-)_+} \mathrm{Fin}_*.$$
*When $m=1$, we omit the adjective "1-subdivided".*

*Remark* 7.2. In Definition 7.1, when $\mathcal{E}$ is a 1-category, an $\mathcal{E}$-valued $m$-subdivided trace theory $T : \Lambda_m^{\mathrm{op}} \times_{\mathrm{Fin}_*} \mathcal{C}^{\otimes} \to \mathcal{E}$ is completely determined by the composite
$$\mathcal{C} \xrightarrow{([1]_{\Lambda_m}, -)} \Lambda_m^{\mathrm{op}} \times_{\mathrm{Fin}_*} \mathcal{C}^{\otimes} \xrightarrow{T} \mathcal{E}$$
along with a natural isomorphism $(\tau_{M,N} : T(M \otimes N) \to T(N \otimes M))_{(M,N) \in \mathcal{C} \times \mathcal{C}}$. Under this description, an $m$-subdivided trace theory descent to a trace theory is a condition (unlike in the $\infty$-categorical setting, it comprises extra data).

Roughly speaking, [NS18, Prop III.3.6 & Lem III.3.7 & Cor III.3.8] (summarized in [KMN23, Thm 6.29]) gives us a $BC_p$-equivariant version of the functor $(-)^{\otimes p} : \mathrm{Sp} \to \mathrm{Sp}^{BC_p}$. We now give a genuine equivariant analogue.

**Notation 7.3.** *Let $G$ be a finite group, and $\mathcal{C}$ a $G$-symmetric monoidal $\infty$-category. Then we will denote $\mathcal{C}_{\mathrm{act}}^{\otimes} := \mathrm{Fin}_G \times_{\mathrm{Span}(\mathrm{Fin}_G)} \mathcal{C}^{\otimes}$.*

*Construction* 7.4. Let $G$ be a finite group, and $\mathcal{C}$ a $G$-symmetric monoidal $\infty$-category. We construct a functor
$$\mathrm{Free}(G) \times_{\mathrm{Fin}} (\mathcal{C}_e)_{\mathrm{act}}^{\otimes} \longrightarrow (\mathcal{C}_G)_{\mathrm{act}}^{\otimes}$$
such that the composite
$$(\mathcal{C}_e)_{\mathrm{act}}^{\otimes} \longrightarrow \mathrm{Free}(G) \times_{\mathrm{Fin}} (\mathcal{C}_e)_{\mathrm{act}}^{\otimes} \longrightarrow (\mathcal{C}_G)_{\mathrm{act}}^{\otimes}$$
is simply given by
$$(S, (X_s)_{s \in S}) \longmapsto (S, (N_e^G X_s)_{s \in S}).$$
Indeed, recall that $(\mathcal{C}_e)_{\mathrm{act}}^{\otimes} = \mathrm{Fin} \times_{\mathrm{free}, \mathrm{Fin}_G} \mathcal{C}_{\mathrm{act}}^{\otimes}$, and $(\mathcal{C}_G)_{\mathrm{act}}^{\otimes} = \mathrm{Fin} \times_{\mathrm{triv}, \mathrm{Fin}_G} \mathcal{C}^{\otimes}$. Consequently, $\mathrm{Free}(G) \times_{\mathrm{Fin}} (\mathcal{C}_e)_{\mathrm{act}}^{\otimes} = \mathrm{Free}(G) \times_{\mathrm{Fin}_G} \mathcal{C}_{\mathrm{act}}^{\otimes}$, and the functor in question is given by the parallel transport of the coCartesian fibration $\mathcal{C}_{\mathrm{act}}^{\otimes} \to \mathrm{Fin}_G$ along the natural transformation from the inclusion $\mathrm{Free}(G) \hookrightarrow \mathrm{Fin}_G$ to the composite functor
$$\mathrm{Free}(G) \xrightarrow{[(-)/G]} \mathrm{Fin} \xrightarrow{\mathrm{triv}} \mathrm{Fin}_G.$$
This natural transformation is pointwisely given by the map $X \to [X/G]_{\mathrm{triv}}$ of finite $G$-sets. By construction, the functor $\mathrm{Free}(G) \times_{\mathrm{Fin}} (\mathcal{C}_e)_{\mathrm{act}}^{\otimes} \to (\mathcal{C}_G)_{\mathrm{act}}^{\otimes}$ carries coCartesian edges over $\mathrm{Free}(G)$ to coCartesian edges over $\mathrm{Fin}$, via $\mathrm{Free}(G) \xrightarrow{[(-)/G]} \mathrm{Fin}$.

*Example* 7.5. Let $G$ be a finite group, and $\mathcal{C}$ an $\infty$-category. Then the construction $[G/H] \mapsto \mathcal{C}^{BH}$ upgrades to a $G$-symmetrc monoidal $\infty$-category, cf. [Hil22, §2]. Construction 7.4 gives rise to a functor
$$\mathrm{Free}(G) \times_{\mathrm{Fin}} \mathcal{C}_{\mathrm{act}}^{\otimes} \longrightarrow (\mathcal{C}^{BG})_{\mathrm{act}}^{\otimes}$$

which carries $(G, X)$ to $X^{\otimes |G|} \in \mathcal{C}^{BG}$ for every $X \in \mathcal{C}$. In fact, this functor is $BG$-equivariant when $G$ is abelian.

Slightly more generally, let $A$ be a commutative ring in $\mathcal{C}^{BG}$, which gives rise to a $G$-$\mathbb{E}_\infty$-algebra in the $G$-symmetric monoidal $\infty$-category $[G/H] \mapsto \mathcal{C}^{BH}$, and thus the construction $[G/H] \mapsto \mathrm{Mod}_A(\mathcal{C}^{BH})$ carries a $G$-symmetric monoidal $\infty$-structure. Then Construction 7.4 gives rise to a functor

$$\mathrm{Free}(G) \times^{\otimes}_{\mathrm{Fin}} \mathrm{Mod}_A(\mathcal{C})^{\otimes}_{\mathrm{act}} \longrightarrow \mathrm{Mod}_A(\mathcal{C}^{BG})^{\otimes}_{\mathrm{act}}$$

which carries $(G, M)$ to $M^{\otimes |G|} \otimes_{A^{\otimes |G|}} A \in \mathrm{Mod}_A(\mathcal{C}^{BG})$ for every $M \in \mathrm{Mod}_A(\mathcal{C})$, where the map $A^{\otimes |G|} \to A$ is informally given by $\bigotimes_{g \in G} a_g \mapsto \prod_{g \in G} g a_g$, cf. [NS18, after Cor IV.2.4]. However, in general, the $\infty$-category $\mathrm{Mod}_A(\mathcal{C}^{BG})^{\otimes}_{\mathrm{act}}$ does not carry a $BG$-action when $G$ is abelian.

*Construction* 7.6. Let $m \in \mathbb{N}_{>0}$ be a positive integer, and $\mathcal{C}$ a $C_m$-symmetric monoidal $\infty$-category. Construction 7.4 along with the commutative diagram

$$\begin{array}{ccc} \Lambda_m^{\mathrm{op}} & \longrightarrow & \Lambda^{\mathrm{op}} \\ \downarrow & & \downarrow \\ \mathrm{Free}(C_m) & \longrightarrow & \mathrm{Fin} \end{array}$$

gives rise to a composite functor

$$\Lambda_m^{\mathrm{op}} \times_{\mathrm{Fin}} (\mathcal{C}_e)^{\otimes}_{\mathrm{act}} \longrightarrow (\mathcal{C}_{C_m})^{\otimes}_{\mathrm{act}} \xrightarrow{\otimes} \mathcal{C}_{C_m}$$

which maps coCartesian edges over $\Lambda_m^{\mathrm{op}}$ to coCartesian edges over Fin, then to equivalences, thus it is an $m$-subdivided trace theory.

*Example* 7.7. Let $r \in \mathbb{N}$, and $m := p^r$. Apply Construction 7.6 to the inclusion $\mathrm{Perm}_{(-)} \subseteq \mathrm{Mod}_{\underline{\mathbb{Z}}}(\mathrm{Sp}^{g(-)})$ of $C_{p^r}$-symmetric monoidal $\infty$-category, we get a $p^r$-subdivided trace theory

$$\Lambda_m^{\mathrm{op}} \times_{\mathrm{Fin}_*} D(\mathbb{Z})^{\otimes} \longrightarrow \mathrm{Mod}_{\underline{\mathbb{Z}}}(\mathrm{Sp}^{gBC_{p^r}})$$

whose restriction to the symmetric monoidal full subcategory $\mathrm{Mod}_{\mathbb{Z}}^{\mathrm{free,fg}} \subseteq D(\mathbb{Z})$ is the $m$-divided trace theory

$$\Lambda_m^{\mathrm{op}} \times_{\mathrm{Fin}_*} (\mathrm{Mod}_{\mathbb{Z}}^{\mathrm{free,fg}})^{\otimes} \longrightarrow \mathrm{Perm}_{BC_{p^r}}$$

which corresponds to the functor

$$\begin{array}{rcl} \mathrm{Mod}_{\mathbb{Z}}^{\mathrm{free,fg}} & \longrightarrow & \mathrm{Perm}_{BC_{p^r}} \\ M & \longmapsto & M^{\otimes p^r} \end{array}$$

along with isomorphisms

$$\begin{array}{rcl} \tau_{M,N} : (M \otimes N)^{\otimes p^r} & = & M \otimes N \otimes \cdots \otimes M \otimes N \\ & = & M \otimes (N \otimes \cdots \otimes M \otimes N) \\ & \xrightarrow{\simeq} & (N \otimes \cdots \otimes M \otimes N) \otimes M \\ & = & N \otimes M \otimes \cdots \otimes N \otimes M \\ & = & (N \otimes M)^{\otimes p^r} \end{array}$$

for every pair $(M, N) \in (\mathrm{Mod}_{\mathbb{Z}}^{\mathrm{free,fg}})^2$ under the description in Remark 7.2, which is denoted by $i^{(r)*} T^{\natural}$ in [Kal18, §4.1].

*Example* 7.8. Let $r \in \mathbb{N}_{>0}$, and $m := p^{r-1}$. Apply Construction 7.6 to the $C_{p^{r-1}}$-symmetric monoidal base change functor

$$\mathrm{Mod}_{\underline{\mathbb{Z}}}(\mathrm{Sp}^{g(-)}) \longrightarrow \mathrm{Mod}_{\underline{W}(\mathbb{F}_p)}(\mathrm{Sp}^{g(-)}),$$

we get two $m$-subdivided trace theories

$$\Lambda_m^{\mathrm{op}} \times_{\mathrm{Fin}_*} D(\mathbb{Z})^{\otimes} \longrightarrow \mathrm{Mod}_{\underline{\mathbb{Z}}}(\mathrm{Sp}^{gBC_{p^{r-1}}})$$
$$\Lambda_m^{\mathrm{op}} \times_{\mathrm{Fin}_*} D(\mathbb{F}_p)^{\otimes} \longrightarrow \mathrm{Mod}_{\underline{W}(\mathbb{F}_p)}(\mathrm{Sp}^{gBC_{p^{r-1}}})$$

which fit into a commutative diagram

$$\begin{array}{ccc} \Lambda_m^{\mathrm{op}} \times_{\mathrm{Fin}_*} D(\mathbb{Z})^{\otimes} & \longrightarrow & \mathrm{Mod}_{\underline{\mathbb{Z}}}(\mathrm{Sp}^{gBC_{p^{r-1}}}) \\ \downarrow & & \downarrow \\ \Lambda_m^{\mathrm{op}} \times_{\mathrm{Fin}_*} D(\mathbb{F}_p)^{\otimes} & \longrightarrow & \mathrm{Mod}_{\underline{W}(\mathbb{F}_p)}(\mathrm{Sp}^{gBC_{p^{r-1}}}) \end{array}.$$

We now upgrade Corollary 6.6 to an equivalence of trace theories.

**Proposition 7.9.** *Let $r \in \mathbb{N}_{>0}$ and $m := p^{r-1}$. The $m$-subdivided trace theory*

$$\Lambda_m^{\mathrm{op}} \times_{\mathrm{Fin}_*} D(\mathbb{F}_p)^{\otimes} \longrightarrow \mathrm{Mod}_{\underline{W}(\mathbb{F}_p)}(\mathrm{Sp}^{gBC_{p^{r-1}}}) \xrightarrow{(-)^{C_{p^{r-1}}}} D(W_r(\mathbb{F}_p)),$$

*where the first functor is the $m$-subdivided trace theory associated to the $C_{p^{r-1}}$-symmetric monoidal $\infty$-category $\mathrm{Mod}_{\underline{W}(\mathbb{F}_p)}(\mathrm{Sp}^{g(-)})$ by Construction 7.6 (cf. Example 7.8), descends canonically to a trace theory*

$$\Lambda^{\mathrm{op}} \times_{\mathrm{Fin}_*} D(\mathbb{F}_p)^{\otimes} \longrightarrow D(W_r(\mathbb{F}_p)),$$

*and its restriction to the symmetric monoidal full subcategory $\mathrm{Mod}_{\mathbb{F}_p}^{\mathrm{free,fg}} \subseteq D(\mathbb{F}_p)$ coincides with Kaledin's trace theory $W_r^{\flat}$ in [Kal18, Prop 4.3].*

*Proof.* The $m$-subdivided trace theory in question preserves fiberwise sifted colimits and is fiberwise polynomial separately in each variable, thus we could simply restrict to the full subcategory $\Lambda_m^{\mathrm{op}} \times_{\mathrm{Fin}_*} \mathrm{Mod}_{\mathbb{F}_p}^{\mathrm{free,fg}}$, that is to say, it is an $m$-subdivided trace theory on the 1-category $\mathrm{Mod}_{\mathbb{F}_p}^{\mathrm{free,fg}}$ of finite free $\mathbb{F}_p$-modules. Moreover, the $m$-subdivided trace theory in question takes values in a 1-category (Corollary 6.6). In this case, an $m$-subdivided trace theory descending to a trace theory is a property (in place of an extra structure), and one can check the coincidence with Kaledin's $W_r^{\flat}$ by explicit computations. □

Now we deduce a consequence for the (Hochschild–)Witt trace theory.

*Remark* 7.10. (**cf. [DKNP23, Ex 2.3]**) Let $r \in \mathbb{N}$. Then the $p^r$-subdivided trace theory obtained by applying Construction 7.6 to the $C_{p^r}$-symmetric monoidal $\infty$-category $\mathrm{Sp}^{g(-)}$ coincides with the $p^r$-subdivided polygonic $\mathrm{THH}(\mathbb{S};-)$. By considering the $C_{p^r}$-symmetric monoidal functor $\mathrm{Sp}^{g(-)} \to \mathrm{Mod}_{\underline{W}(\mathbb{F}_p)}(\mathrm{Sp}^{g(-)})$, we see that the $p^r$-subdivided trace theory

$$\Lambda_m^{\mathrm{op}} \times_{\mathrm{Fin}_*} D(\mathbb{F}_p)^{\otimes} \longrightarrow \mathrm{Mod}_{\underline{W}(\mathbb{F}_p)}(\mathrm{Sp}^{g(-)})$$

coming from the $C_{p^r}$-symmetric monoidal $\infty$-category $\mathrm{Mod}_{\underline{W}(\mathbb{F}_p)}(\mathrm{Sp}^{g(-)})$ can be identified with the the $p^r$-subdivided trace theory

$$\mathrm{THH}(\mathbb{S};-) \otimes_{\mathrm{THH}(\mathbb{S};\mathbb{F}_p)}^{\mathbb{L}} \underline{W}(\mathbb{F}_p) \simeq \mathrm{THH}(\mathbb{F}_p;-) \otimes_{\mathrm{THH}(\mathbb{F}_p)}^{\mathbb{L}} \underline{W}(\mathbb{F}_p)$$

(this equivalence follows from the parametrized symmetric monoidal structure on the trace theory THH) on $D(\mathbb{F}_p)$. Arguing as in Proposition 7.9 (namely, restricting to the full subcategory $\mathrm{Mod}_{\mathbb{F}_p}^{\mathrm{free,fg}} \subseteq D(\mathbb{F}_p)$, and using the polynomiality of functors), we see that the $p^r$-subdivided trace theory $(\mathrm{THH}(\mathbb{F}_p;-) \otimes_{\mathrm{THH}(\mathbb{F}_p)}^{\mathbb{L}} \underline{W}(\mathbb{F}_p))^{C_{p^r}}$ descends to a trace theory, which coincides with Kaledin's trace theory $W_{r+1}^{\natural}$ as well.

Then it follows from the trace-theory formalism that

**Proposition 7.11.** *Let $A$ be an $\mathbb{E}_1$-$\mathbb{F}_p$-algebra, and $r \in \mathbb{N}_{>0}$. Then the $r$-truncated Hochschild–Witt homology $W_r \mathrm{HH}(A)$ of $A$ (defined in [Kal19, Def 4.1]) can be identified with the $r$-truncated Hochschild–Witt homology*

$$W_r \mathrm{HH}(A/\mathbb{F}_p) \in D(W_r(\mathbb{F}_p))^{B(\mathbb{T}/C_{p^{r-1}})}$$

*of $A$ relative to $\mathbb{F}_p$.*

To complete the comparison, in view of Corollary 4.2 and Proposition 7.11, it suffices to identify the map $R: W_{r+1} \mathrm{HH}(-) \to W_r \mathrm{HH}(-)$ with the map

$$(\mathrm{THH}(-) \otimes_{\mathrm{THH}(\mathbb{F}_p)}^{\mathbb{L}} \underline{W}(\mathbb{F}_p))^{C_{p^r}} \longrightarrow (\mathrm{THH}(-) \otimes_{\mathrm{THH}(\mathbb{F}_p)}^{\mathbb{L}} \underline{W}(\mathbb{F}_p))^{C_{p^{r-1}}}$$

out of the pre-cyclotomic structure on $\mathrm{THH}(-) \otimes_{\mathrm{THH}(\mathbb{F}_p)}^{\mathbb{L}} \underline{W}(\mathbb{F}_p)$. As for Proposition 7.11, we identify the map $W_{r+1}^{\natural} \to W_r^{\natural}$ of trace theories with the map

$$(\mathrm{THH}(\mathbb{F}_p;-) \otimes_{\mathrm{THH}(\mathbb{F}_p)}^{\mathbb{L}} \underline{W}(\mathbb{F}_p))^{C_{p^r}} \longrightarrow (\mathrm{THH}(\mathbb{F}_p;-) \otimes_{\mathrm{THH}(\mathbb{F}_p)}^{\mathbb{L}} \underline{W}(\mathbb{F}_p))^{C_{p^{r-1}}}$$

of trace theories, on $\mathrm{Mod}_{\mathbb{F}_p}^{\mathrm{free,fg}}$, and thus it suffice to identify this map with the map

$$(\mathrm{THH}(\mathbb{F}_p;-) \otimes_{\mathrm{THH}(\mathbb{F}_p)}^{\mathbb{L}} \underline{W}(\mathbb{F}_p))^{C_{p^r}} \longrightarrow (\mathrm{THH}(\mathbb{F}_p;-) \otimes_{\mathrm{THH}(\mathbb{F}_p)}^{\mathbb{L}} \underline{W}(\mathbb{F}_p))^{C_{p^{r-1}}} \quad (7.1)$$

coming from the pre-polygonic structure on $\mathrm{THH}(\mathbb{F}_p;-) \otimes_{\mathrm{THH}(\mathbb{F}_p)}^{\mathbb{L}} \underline{W}(\mathbb{F}_p)$ as natural transformations of functors out of $\mathrm{Mod}_{\mathbb{F}_p}^{\mathrm{free,fg}}$.

By previous computations, on $\mathrm{Mod}_{\mathbb{F}_p}^{\mathrm{free,fg}}$, we can replace the map (7.1) by applying $\pi_0$ to it, and thus we can identify it with taking $\pi_0$ to the map

$$\mathrm{THH}(\mathbb{F}_p;-)^{C_{p^r}} \longrightarrow \mathrm{THH}(\mathbb{F}_p;-)^{C_{p^{r-1}}}$$

coming from the polygonic structure on $\mathrm{THH}(\mathbb{F}_p;-)$, which is subsequently identified with $W_{r+1}^{\natural} \to W_r^{\natural}$ in [DKNP23].

In summary, we have

**Proposition 7.12.** *Let $A$ be an $\mathbb{E}_1$-$\mathbb{F}_p$-algebra. Then the sequential system of Hochschild–Witt homology $(W_r \mathrm{HH}(A))_{r \in \mathbb{N}_{>0}}$ of $A$ (defined in [Kal19, Def 4.1]), with restriction maps are transition maps, can be identified with Hochschild–Witt homology*

$$\left(W_r \mathrm{HH}(A/\mathbb{F}_p) \in D(W_r(\mathbb{F}_p))^{B(\mathbb{T}/C_{p^{r-1}})}\right)_{r \in \mathbb{N}_{>0}}$$

*relative to $\mathbb{F}_p$, with transition maps coming from the pre-$p$-cyclotomic structure.*

*Remark* 7.13. The identification above of maps $R$ is dirty. In the future, we will give another argument by seeking a suitable generalization without the identification of $\pi_0 \mathrm{THH}^{C_{p^{r-1}}}$ in [DKNP23].

Appendix A. $G$-typical Witt vectors as the relative norm

In this appendix, we briefly show that Thomas Read's $G$-typical Witt vectors with coefficients [Rea23] can be identified with the relative norm under flatness, as a simple consequence of base-change formulae. For sake of simplicity, we only consider finite groups $G$, and without truncation sets. We start with the case that the commutative ring is $\mathbb{Z}$.

*Construction* A.1. Let $M$ be a connective $\mathbb{Z}$-module spectrum. Then the 0th Postnikov truncation map $N_e^G M \to \tau_{\leqslant 0}(N_e^G M)$ in $\mathrm{Mod}_{N_e^G \mathbb{Z}}(\mathrm{Sp}^{gBG})$ induces a map

$$N_e^G M \otimes_{N_e^G \mathbb{Z}}^{\mathbb{L}} \tau_{\leqslant 0}(N_e^G \mathbb{Z}) \longrightarrow \tau_{\leqslant 0}(N_e^G M)$$

in $\mathrm{Mod}_{\pi_0(N_e^G \mathbb{Z})}(\mathrm{Sp}^{gBG})$ which is functorial in $M \in D(\mathbb{Z})_{\geqslant 0}$, and becomes an equivalence after taking $\tau_{\leqslant 0}$.

**Lemma A.2.** *Let $M$ be a flat $\mathbb{Z}$-module. Then the map in Construction A.1 becomes an equivalence.*

*Proof.* It suffices to see that the $G$-spectrum $N_e^G M \otimes_{N_e^G \mathbb{Z}}^{\mathbb{L}} \tau_{\leqslant 0}(N_e^G \mathbb{Z})$ is concentrated in degree 0. Since every flat module is a filtered colimit of finite free modules (Lazard's theorem), without loss of generality, we may assume that the $\mathbb{Z}$-module $M$ is finite free. We write $M = \mathbb{Z}[S]$ for a finite set $S$. Then

$$\begin{aligned} N_e^G M \otimes_{N_e^G \mathbb{Z}}^{\mathbb{L}} \tau_{\leqslant 0}(N_e^G \mathbb{Z}) &\simeq (S^{\times G} \otimes N_e^G \mathbb{Z}) \otimes_{N_e^G \mathbb{Z}}^{\mathbb{L}} \tau_{\leq 0}(N_e^G \mathbb{Z}) \\ &\simeq S^{\times G} \otimes \tau_{\leq 0}(N_e^G \mathbb{Z}) \end{aligned}$$

is concentrated in degree 0, where $S^{\times G}$ is viewed as a finite $G$-set, and $S^{\times G} \otimes N_e^G \mathbb{Z}$ is the $\mathrm{Fin}_G$-action on $\mathrm{Sp}^{gBG}$. $\square$

It follows from [Rea23, Thm A] that

**Corollary A.3.** *Let $M$ be a flat $\mathbb{Z}$-module. Then there is an equivalence*

$$N_e^G M \otimes_{N_e^G \mathbb{Z}}^{\mathbb{L}} \underline{W}_G(\mathbb{Z}) \xrightarrow{\simeq} \underline{W}_G(\mathbb{Z}; M)$$

*in $\mathrm{Mod}_{\underline{W}_G(\mathbb{Z})}(\mathrm{Sp}^{gBG})$, which is functorial in $M \in \mathrm{Mod}_{\mathbb{Z}}^{\flat}$.*

Now we deduce the general result.

**Proposition A.4.** *Let $R$ be a commutative ring, and $M$ a flat $R$-module. Then there is an equivalence*

$$N_e^G M \otimes_{N_e^G R}^{\mathbb{L}} \underline{W}_G(R) \xrightarrow{\simeq} \underline{W}_G(R; M) \tag{A.1}$$

*in $\mathrm{Mod}_{\underline{W}_G(R)}(\mathrm{Sp}^{gBG})$, which is functorial in $(R, M) \in \mathrm{Mod}^{\flat}$.*

*Proof.* By Corollary A.3, we can rewrite the map as

$$\underline{W}_G(\mathbb{Z}; M) \otimes_{\underline{W}_G(\mathbb{Z}; R)}^{\mathbb{L}} \underline{W}_G(R) \longrightarrow \underline{W}_G(R; M).$$

By [Rea23, Cor 4.21], this map becomes an equivalence after taking $\tau_{\leqslant 0}$. It suffices to show that the source of (A.1) is concentrated in degree 0. Since every $R$-module is a filtered colimit of finite free $R$-modules, without loss of generality, we may assume that $M$ is finite free. We write $R = \mathbb{Z}[S]$ for a finite set $G$. Then

$$\begin{aligned} N_e^G M \otimes_{N_e^G R}^{\mathbb{L}} \underline{W}_G(R) &\simeq S^{\times G} \otimes N_e^G R \otimes_{N_e^G R}^{\mathbb{L}} \underline{W}_G(R) \\ &\simeq S^{\times G} \otimes \underline{W}_G(R) \end{aligned}$$

is concentrated in degree 0 since so is $\underline{W}_G(R)$. □